\documentclass[10pt]{amsart}

\usepackage[all]{xypic}

\usepackage{latexsym}

\usepackage{amssymb}
\usepackage{amsfonts}
\usepackage{amscd}
\usepackage{amsmath,amsthm}

\newtheorem{lemma}{Lemma}[section]

\newtheorem{lem}[lemma]{Lemma}
\newtheorem{prop}[lemma]{Proposition}
\newtheorem{thm}[lemma]{Theorem}
\newtheorem{cor}[lemma]{Corollary}

{

}

\theoremstyle{definition}

\theoremstyle{remark}

\newtheorem{remark}[lemma]{Remark}

\newcommand{\appsectionA}[1]{\let\oldthesection\thesection
   \renewcommand{\thesection}{Appendix A}
  \section{#1}\let\thesection\oldthesection}

\newcommand{\appsectionB}[1]{\let\oldthesection\thesection
   \renewcommand{\thesection}{Appendix B}
  \section{#1}\let\thesection\oldthesection}

\numberwithin{equation}{section}

\newenvironment{pf}{\noindent{\bf Proof.}}{\hfill $\square$\medskip}

\def\CC{{\mathbb C}}

\def\FF{{\mathbb F}}

\def\NN{{\mathbb N}}

\def\PP{{\mathbb P}}

\def\RR{{\mathbb R}}

\def\ZZ{{\mathbb Z}}

\def\Sol{{\bar S}}

\def\0ol{{\bar 0}}
\def\1ol{{\bar 1}}
\def\2ol{{\bar 2}}
\def\ol2{{\bar 2}}
\def\3ol{{\bar 3}}
\def\4ol{{\bar 4}}
\def\5ol{{\bar 5}}
\def\6ol{{\bar 6}}
\def\7ol{{\bar 7}}
\def\8ol{{\bar 8}}
\def\9ol{{\bar 9}}

\def\bold0{{\bf 0}}
\def\bold1{{\bf 1}}
\def\bold2{{\bf 2}} 
\def\bold3{{\bf  3}}
\def\bold4{{\bf 4}}
\def\bold5{{\bf 5}}
\def\bold6{{\bf 6}}
\def\bold7{{\bf 7}}
\def\bold8{{\bf 8}}
\def\bold9{{\bf 9}}

\def\vul{{\underline{v}}}

\def\P2Skly{\PP^2_{Skly}}

\def\coker{\operatorname {coker}}

\def\diag{\operatorname {diag}}
\def\End{\operatorname {End}}
\def\Ext{\operatorname {Ext}}

\def\gr{\operatorname {gr}}

\def\Hom{\operatorname {Hom}}

\def\ker{\operatorname {ker}}

\def\th{\operatorname {th}}    

\def\dim{\operatorname{dim}}

\def\End{\operatorname{End}}

\def\Ext{\operatorname{Ext}}
\def\Fdim{{\sf Fdim}}
\def\fdim{{\sf fdim}}

\def\Gr{{\sf Gr}}

\def\gr{{\sf gr}}

\def\Hom{\operatorname{Hom}}

\def\limproj{\varprojlim}
\def\liminj{\varinjlim}

\def\min{\operatorname{min}}
\def\mod{{\sf mod}}
\def\Mod{{\sf Mod}}

\def\Proj{\operatorname{Proj}}
\def\Projnc{\operatorname{Proj}_{nc}}

\def\QGr{\operatorname{\sf QGr}}
\def\qgr{\operatorname{\sf qgr}}

\def\ul1{\operatorname{\underline{1}}}

\def\l{\leftarrow}

\def\a{\alpha}
\def\b{\beta}
\def\c{\gamma}

\def\l{\lambda}

\def\s{\sigma}

\def\ve{\varepsilon}

\def\L{\Lambda}

\def\fm{{\mathfrak m}}

\def\sA{{\sf A}}

\def\sF{{\sf F}}

\def\sP{{\sf P}}

\def\sT{{\sf T}}

\def\sX{{\sf X}}

\def\cal{\mathcal}

\def\cF{{\cal F}}
\def\cG{{\cal G}}

\def\cM{{\cal M}}
\def\cN{{\cal N}}
\def\cO{{\cal O}}

\def\cT{{\cal T}}

\def\coh{{\sf coh}}

\def\Qcoh{{\sf Qcoh}}

\def\dirlim{\mathop{\vtop{\baselineskip -100pt\lineskip -1pt\lineskiplimit 0pt
\setbox0\hbox{lim}\copy0\hbox to \wd0{\rightarrowfill}}}\limits}
\def\invlim{\mathop{\vtop{\baselineskip -100pt\lineskip -1pt\lineskiplimit 0pt
\setbox0\hbox{lim}\copy0\hbox to \wd0{\leftarrowfill}}}\limits}

\def\I11{{1 \kern -0.8pt \! \mbox{l}}}
\def\mumu{{\mu\kern-4.2pt\mu}}
\def\bfmu{{\mu\kern-4.2pt\mu}}
\def\2slash{\backslash \! \backslash}

\def\boxtimes{\setbox0\hbox{$\Box$}\copy0\kern-\wd0\hbox{$\times$}}

\def\ldot{{\:\raisebox{2pt}{\text{\circle*{1.5}}}}}

\begin{document}

\title[The space of Penrose tilings and $k\langle x,y\rangle/(y^2)$]{The space of Penrose tilings and the non-commutative curve with homogeneous coordinate ring $k\langle x,y\rangle/(y^2)$}

\author{S. Paul Smith}
\address{ Department of Mathematics, Box 354350, Univ.  Washington, Seattle, WA 98195}
\email{smith@math.washington.edu}

\subjclass{14A22, 16E20, 16E50, 16S38, 16W50, 19K14, 52C20, 52C22, 58B34}

\keywords{Homogeneous coordinate ring, non-commutative algebraic geometry, Penrose tiling, 
aperiodic tilings, AF-algebra, Grothendieck group,}

\thanks{The author was partially supported by the 
National Science Foundation, Award No. 0602347}


\begin{abstract}
The non-commutative algebraic curve with homogeneous coordinate ring $\CC\langle x,y\rangle/(y^2)$ 
is shown to be a non-commutative algebraic-geometric analogue of the space of Penrose tilings of the plane. 
Individual tilings determine ``points'' on the non-commutative curve and the tilings coincide under isometry if and only if the skyscaper sheaves of the corresponding points are isomorphic. The category of quasi-coherent sheaves on the curve is equivalent to the category of modules over a von Neumann regular ring that is a direct limit of finite dimensional semisimple algebras. The norm closure of this von Neumann regular ring is the AF-algebra that Connes associates to the space of Penrose tilings. There is an algebraic analogue of the fact that every isometry-invariant subset of tilings is dense in the set of all Penrose tilings.
\end{abstract}

\maketitle

\section{Introduction}

Let $k$ be a field and $k\langle x,y\rangle$ the free algebra on $x$ and $y$.

Let $B:=k\langle x,y\rangle/(y^2)$ have the $\ZZ$-grading $\deg x =\deg y=1$. 

The Fibonacci sequence $(f_n)_n$ is defined by $f_0=f_1=1$ and $f_{n+1}:=f_n+f_{n-1}$ for $n \ge 1$.

\subsection{}
This paper concerns the algebras $k\langle x,y\rangle/(y^{r+1})$ for $r \ge 1$. For simplicity we will only discuss the case $r=1$ in this introduction, that being the case related to Penrose tilings. However, all
results in this introduction extend to $r \ge 2$. When $r  \ge 2$ the space of Penrose tilings can be replaced either by (a) the space of (0,1)-sequences having at most $r$ consecutive 1s modulo the equivalence relation that two such sequences are equivalent if they are eventually the same or by 
(b) an associated space of aperiodic tilings in $\RR^r$. Many papers discuss the connection between (a) and (b) but \cite{BJS} is a good reference for those interested in non-commutative algebraic geometry because it emphasizes the relationship with infinite paths in Bratteli diagrams and 
equality between the tails of such paths is related to equivalence of tilings under isometry and to isomorphisms between the images of point modules in the quotient category $\QGr B$ (this terminology and notation is explained in sections \ref{intro.main} and \ref{ssect.2nd.thm}). 

\subsection{}
Penrose tilings are defined in Appendix A. The appendix gives a precise description
of the relation between Penrose tilings and (0,1)-sequences.

\subsection{}
\label{sect.1.1}
The space of Penrose tilings of the plane is an instructive example showing that  
non-commutative algebras can capture important features of geometric objects that are not susceptible to the usual topological, analytic, and differential geometric, tools.   A nice account of this example appears in  Connes book \cite[Sect. II.3]{NCG}. Using a  general principle that we elaborate on in section \ref{sect.inv.lim.dlim}, Connes associates to the space of Penrose tilings  the C${}^*$-algebra 
$$
\Sol:= \liminj_n (S_n,\phi_n), 
$$ 
the  direct limit in the category of C${}^*$-algebras of the finite dimensional C${}^*$-algebras  
 $$
 S_n := M_{f_n}(\CC) \oplus M_{f_{n-1}}(\CC)
 $$
where the maps  $\phi_n:S_n \to S_{n+1}$ in the directed system are  
$$
\phi_n(a,b) := \Bigg(\begin{pmatrix}
   a   &  0  \\
   0   &  b
\end{pmatrix}, a \Bigg).
$$

\subsection{The main result}
\label{intro.main}
Non-commutative {\it algebraic} geometry concerns abelian categories that behave like categories of quasi-coherent sheaves on algebraic varieties. We view $B$ as a homogeneous coordinate ring 
of a ``non-commutative scheme'', $\Projnc B$. By fiat, the category of ``quasi-coherent sheaves'' on $\Projnc B$  is   
$$
\QGr B:= \frac{\Gr B}{\Fdim B}\, ,
$$
the quotient of the category $\Gr B$ of $\ZZ$-graded left $B$-modules modulo the Serre subcategory of modules that are unions of their finite dimensional submodules.

By definition, the objects in $\QGr B$ are the same as those in $\Gr B$ but there are more isomorphisms in $\QGr B$. For example, every finite dimensional graded $B$-module is isomorphic to the zero object in $\QGr B$.

We write $\cO$ for $B$ as an object in $\QGr B$. It is the ``structure sheaf'' 
for    $\Projnc B$. (If $R$ is the commutative polynomial ring on $n$ variables with its standard grading, then $\QGr R$ is equivalent to $\Qcoh \PP^{n-1}$, the category of quasi-coherent sheaves on projective $(n-1)$-space, and $R$ does correspond to the  structure sheaf $\cO_{\PP^{n-1}}$.)

To avoid confusion about which category an object belongs to we introduce notation for the canonical quotient functor, namely
$$
\b^*:\Gr B \to \QGr B.
$$
Thus $\cO=\b^*B$.

Since $B$ is a coherent ring (Proposition \ref{prop.B.coh}), $\Gr B$ and $\QGr B$ are locally coherent categories. Let $\gr B$ be the full subcategory of $\Gr B$ consisting of finitely presented modules.
It is an abelian subcategory and is closed under extensions in $\Gr B$.  Let $\fdim B$ be the full subcategory of $\Gr B$ consisting of  finite dimensional modules. A finite dimensional graded module is finitely presented because all its composition factors are isomorphic to the finitely presented module
$B/B_{\ge 1}$. Thus $\fdim B = \gr B \cap \Fdim B$. We write 
$$
\qgr B:= \frac{\gr B}{\fdim B}
$$
and note that $\qgr B$ is equivalent to the full subcategory of $\QGr B$ consisting of finitely presented, in the categorical sense, objects.  We think of $\qgr B$ as the category of ``coherent sheaves'' on $\Projnc B$.

\begin{thm}
\label{thm.main}
Let $S$ be the direct limit in the category of $\CC$-algebras of the directed system $(S_n,\phi_n)$
in section \ref{sect.1.1}. (The C${}^*$-algebra $\Sol$ is the norm closure of $S$.) Let $\Mod S$ be the category of right $S$-modules and $\mod S$ its full subcategory of finitely presented right $S$-modules.
There is an equivalence of categories
$$
\QGr \Bigg(\frac{\CC\langle x,y\rangle}{(y^2)}\Bigg) \, \equiv \, \Mod S
$$
that restricts to an equivalence between $\qgr\Big(\CC\langle x,y\rangle/(y^2)\Big)$ and $\mod S$. Under the equivalence $\cO$ corresponds to $S$.
\end{thm}

The ring $S$ is a simple von Neumann regular algebra.
It is left (and right) coherent because it is a direct limit of coherent rings and the homomorphisms in the directed system make each ring flat as a right module over the earlier rings in the directed system. 
Thus, independently of the equivalence of categories, 
one knows that $\Mod S$ is a locally coherent category and 
$\mod S$ is an abelian category.

The field $\CC$ is of no importance in proving Theorem \ref{thm.main}. There is a version of Theorem
\ref{thm.main} for every field $k$. Given Proposition \ref{prop.F2}, it could be argued that 
 $k=\FF_2$ is the most natural choice when considering the relation between $\Projnc B$ and the space of Penrose tilings.

\begin{cor}
\label{cor.inj.proj}
Every object in $\qgr B$ is injective and projective.
\end{cor}

The proof is straightforward: all modules over a von Neumann regular ring are flat \cite[Cor. 1.13]{vNR}
so finitely presented $S$-modules are projective. But $\mod S$ is equivalent to $\qgr B$ so every object
in $\qgr B$ is projective and every short exact sequence in $\qgr B$ splits. 

Because $S$ is a von Neumann regular ring of countable dimension every left ideal in it is countably generated and therefore projective  \cite[Cor. 2.15]{vNR}.  This implies that all subobjects of $\cO$,
not just the finitely presented ones, are projective.

\subsection{Point modules, ``points'' in $\Projnc B$, and Penrose tilings}
\label{ssect.2nd.thm}

As in algebraic geometry, certain simple objects in $\QGr B$ are thought of as
 ``skyscraper sheaves'' at ``points'' of $\Projnc B$. 

In \cite{ATV1}, a 
graded $B$-module $M$ is called a {\sf point module}  if $\dim_k M_n=1$ for all 
$n \ge 0$ and $M=BM_0$. Since $B$ is generated as a $k$-algebra by $B_1$,
every non-zero submodule of a point module $M$ has finite codimension so  $M$ is simple as an
object in $\QGr B$.  

A {\sf Penrose sequence} is a sequence $z=z_0z_1\ldots$ of $0$s and $1$s having no consecutive 1s.
We write ${\bf P}$ for the set of Penrose sequences.
Each $z \in {\bf P}$ determines a Penrose tiling $\cT_z$ and every tiling is obtained in this way
(see Gr\"unbaum and Shepard \cite[10.5.9]{GS} and Appendix A below). Penrose sequences $z$ and $z'$ are {\sf equivalent}, denoted $z \sim z'$,  if there is an $m$ such that $z_n=z'_n$ for all $n \ge m$. We denote this by $z_{\ge m}=z'_{\ge m}$. Given a Penrose sequence $z$ we define the 
point module 
$$
M_z:=ke_0 \oplus ke_1 \oplus \cdots, \qquad \deg e_i=i,
$$
by
$$
x.e_i=(1-z_i) e_{i+1} 
\qquad \hbox{and} \qquad 
y.e_i=z_ie_{i+1}.
$$
Since $y^2.e_i=z_iz_{i+1}e_{i+2} =0$ because $z_iz_{i+1} \ne 11$,  $M_z$ is a graded $B$-module. 
We write $\cO_z$ for $M_z$ viewed as an object in $\QGr B$; i.e., $\cO_z = \b^*M_z$.

\begin{thm}
\label{thm.2}
If $z$ is a Penrose sequence, then $\cO_z$ is a simple object in $\QGr B$ and the following
are equivalent:
\begin{enumerate}
  \item 
  $\cO_z \cong \cO_{z'}$;
  \item 
  $z \sim z'$;
  \item 
  there is an isometry $\s$ of $\RR^2$ such that $\s(\cT_z)=\cT_{z'}$.
\end{enumerate}
\end{thm}

The equivalence between (1) and (2) is routine: point modules $M$ and $M'$ are isomorphic as objects in $\QGr B$ if and only if $M_{\ge n} \cong M'_{\ge n}$ for some $n$; but $(M_z)_{\ge n} \cong (M_{z'})_{\ge n}$ if and only if $z_{\ge n} = z'_{\ge n}$. The equivalence between (2) and (3) can be found in 
\cite[10.5.9]{GS}, for example.

\subsubsection{Degree shift and Serre twist}
If $M$ is a graded $B$-module and $n \in \ZZ$ we define the graded module $M(n)$ to be 
$M$ as a $B$-module but with the grading $M(n)_i:=M_{n+i}$. 
Each $(n)$ is an automorphism of $\Gr B$ and descends to an automorphism of $\QGr B$ that we denote by $\cF \rightsquigarrow \cF(n)$ and call a {\sf Serre twist}.

There is an automorphism $\s$ of ${\bf P}$ defined by $\s(z_0z_1z_2 \ldots):=z_1z_2\ldots$.
It is clear that $M_z(1)_{\ge 0} = M_{\s(z)}$, whence 
$$
\cO_z(1)=\cO_{\s(z)}.
$$
The unique $\cO_z$ such that $\cO_z \cong \cO_z(1)$ corresponds to the cart-wheel
tiling \cite[10.5.14]{GS}. This $\cO_z$ is the structure sheaf of the classical point $\Proj B/(y)$
on $\Projnc B$.

\subsection{Topology of ${\bf P}$ and algebraic properties of $\QGr B$}
\label{ssect.top.P}

We establish a nice parallel between properties of $\Projnc B$ and the topological structure of ${\bf P}$
and ${\bf P}/\!\! \sim$.

The usual topology on ${\bf P}$ is the subspace topology inherited from its being a subset of 
$\{0,1\}^\NN$ endowed with the product topology. It is also the topology induced by the metric
$$
d(z,z')=\sum_{n=0}^\infty 2^{-n}|z_n-z_n'|.
$$

Motivated by Theorems \ref{thm.main} and \ref{thm.2}, and the common view that certain C${}^*$-algebras that are direct limits of semisimple $\CC$-algebras are the ``function rings'' for 
appropriate topological spaces of tilings,  we prove several results that support the view 
that $\Projnc B$ is a non-commutative 
algebraic-geometric version of the space of Penrose tilings.

\subsubsection{Dimension}
Since ${\bf P}/\!\!\sim$ has the coarse topology its topological dimension is zero. This might be seen as corresponding to the fact that every short exact sequence in $\qgr B$ splits (Corollary \ref{cor.inj.proj}).
(A reduced scheme $X$ has this  property if and only if it has dimension 0.)
 
 \subsubsection{Density}
A striking feature of Penrose tilings is that every finite region in one tiling appears infinitely often in {\it every} tiling. A consequence of this is that every equivalence class of tilings, more precisely, every equivalence class of Penrose sequences, is dense in ${\bf P}$. In section \ref{ssect.dense} we explain why the following result is an algebraic analogue of this density result.

\begin{prop}
\label{prop.links}
If $z,z' \in {\bf P}$, then $\Ext^1_{\QGr B}(\cO_z,\cO_{z'}) \ne 0$.
\end{prop}

Proposition \ref{prop.links} follows from Theorem \ref{thm.main} and Herbera's result
that $\Ext^1_R(U,V) \ne 0$ for all simple modules $U$ and $V$ over a simple non-artinian von Neumann regular algebra $R$ having countable dimension over its center \cite[Prop. 2.7]{DH}.

\subsubsection{Closed subspaces}
Because every equivalence class in ${\bf P}$ is dense in ${\bf P}$, the only closed sets in
${\bf P}/\!\!\sim$ with respect to the quotient topology are itself and the empty set. 

There is a notion of a closed subspace in non-commutative algebraic geometry \cite[p.20]{VdB}, \cite[Defn. 2.4]{Sm0}. When $X$ is a noetherian scheme having an ample line bundle (for example, 
when $X$ is quasi-projective over the spectrum of a field) closed subspaces of $X$ in the non-commutative sense are the same things as closed subschemes in the usual sense of algebraic geometry \cite[Thm. 4.1]{Sm0}. Rosenberg has shown that the closed subspaces of an affine non-commutative scheme, i.e., one whose category of quasi-coherent sheaves is equivalent to a module category, say $\Mod R$, are in natural bijection with the two-sided ideals of $R$. Since 
the ring $S$ in Theorem \ref{thm.main} is simple the following result holds.

\begin{prop}
The only closed subspaces of $\Projnc B$ are itself and the empty set.
\end{prop}
 
\subsection{$\Projnc B$ is a ``smooth curve''}

Although $B$ appears in the classification of 2-generator algebras subject to one quadratic relation it 
has not been studied much because it has infinite global homological dimension, is homologically nasty in other ways, is not noetherian, and has exponential rather than polynomial growth. In short, it is 
most unlike the rings that appear in classical algebraic geometry

The ``nice'' 2-generator algebras subject to one quadratic relation are 
$R_J:=k\langle x,y\rangle/(xy-yx+y^2)$ and   
$R_q:=k\langle x,y\rangle/(qxy-yx)$ where $q \in k^\times$.
The categories $\QGr( R_J)$ and $\QGr( R_q)$ are equivalent to $\Qcoh \PP^1$.

Because the ring $S$ in Theorem \ref{thm.main} is von Neumann regular it has global homological dimension one, and every subobject of $\cO$ is projective, both of which suggest that $\Projnc B$ is some kind of ``non-commutative smooth curve''.

\subsection{The Grothendieck group}

We compute $K_0(\qgr(B))$ as an ordered abelian group in section \ref{sect.K0}. 
First, we compute it as an abstract group without using the equivalence with $\mod S$ 
(Theorem \ref{thm.K0.1}). The first step towards doing this is to show that $\qgr(B)$ is equivalent to 
$\qgr(kQ)$ where $kQ$ is the path algebra of the quiver
$$
 \UseComputerModernTips
\xymatrix{
\bullet    \ar@/^/[rr]  \ar@(ul,dl)[]  && \bullet \ar@/^/[ll] 
}
$$
Using the ideas in \cite{HS1} and \cite{HS2} we show there is an injective degree-preserving algebra homomorphism $f:B \to kQ$ with the property that the functor $kQ \otimes_B -$ induces an equivalence of categories $\qgr(B) \stackrel{\equiv}{\longrightarrow}  \qgr(kQ)$ and hence an isomorphism
$K_0(\qgr(B)) \stackrel{\sim}{\longrightarrow} K_0(\qgr(kQ))$.

The advantage of this is that $kQ$ has finite global dimension (=1) so one can compute $K_0(\gr(kQ))$
easily and then use the localization sequence $K_0(\fdim(kQ)) \to K_0(\gr(kQ)) \to K_0(\qgr(kQ)) \to 0$
to compute $K_0(\qgr(kQ))$. The result is that   
$$
K_0(\qgr B) \cong \frac{\ZZ[t,t^{-1}]}{(1-t-t^2)}
$$
with $[\cO(n)] \leftrightarrow t^{-n}$.

The ring homomorphism $K_0(\qgr(B))\to \RR$ given by sending $t^{-1}$ to  the Golden Ratio,
 $\tau=\frac{1+\sqrt{5}}{2}$, is an isomorphism from 
$K_0(\qgr B)$ onto $\ZZ+\ZZ\tau$.  It  sends  $[\cO]$ to 1 and $[\cO(1)]$ to $\tau$.
Theorem \ref{thm.K0} shows that  this isomorphism sends the positive cone in $K_0(\qgr B)$ 
to $\RR^+ \cap(\ZZ + \ZZ \tau)$.  

We use the work of Bratteli \cite{B}  and Elliott \cite{E} on AF-algebras and refer to 
Effros's lecture notes \cite{Eff} for all matters related to AF-algebras.

\subsubsection{The Grothendieck group for $\qgr k\langle x,y\rangle/(y^{r+1})$}

The case $r>1$ is more interesting. The Bratteli diagram for the directed system $(S_n,\phi_n)$ and 
$K_0(S_n,\phi_n)$ is constant, or stationary
as is usually said (its shape is shown in Proposition \ref{prop.Brat.Sn}). This simplifies the computations. The characteristic  polynomial of the matrix (\ref{M.matrix}) for this stationary system 
is $t^{r+1}-t^r-\cdots -t-1$. This polynomial has a unique positive
real root of multiplicity one, $\a$ say,  that is strictly larger than the absolute values of all other 
roots. We show that $K_0(\mod S)$, which is isomorphic to $K_0(\qgr B)$, 
is isomorphic to $\ZZ[\a]$. Under this isomorphism the positive cone of $K_0(\qgr B)$ is 
$\ZZ[\a]  \cap \RR^+$, $[\cO] \leftrightarrow 1$, $[\cO(-n) ]  \leftrightarrow \a^n$, and the twist $[\cF] \rightsquigarrow [\cF(-1)]$ corresponds to multiplication by $\a$.

 \subsection{Acknowledgements}
 It is a pleasure to thank Ken Goodearl, Doug Lind, and Boris Solomyak, for many helpful conversations. I am particularly grateful to Gautam Sisodia for finding a fatal flaw in 
 an earlier ``proof'' that $K_0(\qgr B) \cong \ZZ[\a]$; he also suggested using the results in 
 \cite{HS1} to repair the proof. Sections \ref{sect.kQ.0}, \ref{sect.kQ.1}, and \ref{sect.kQ.2}, 
 were prompted by his comments. 
 
 I thank the US National Science Foundation for partial support under NSF grant 0602347.

 \section{Some general results}
 
 \subsection{Graded-coherent rings}
 
A ring is {\sf left coherent} if all its finitely generated left ideals are finitely presented.
There is a similar notion on the right.

In this section $A$ and $B$ are arbitrary $\ZZ$-graded $k$-algebras. We write $\Gr A$ for the category of graded left $A$-modules. We
say that $A$ is {\sf  left graded-coherent} if every finitely generated graded left ideal of $A$ is 
finitely presented. If $A$ is left graded-coherent we write $\gr A$ for the full subcategory of $\Gr A$ consisting of finitely presented modules. The category $\gr A$ is an abelian subcategory of $\Gr A$ and is closed 
under extensions. If $M \in \gr A$ and $L$ is a finitely generated graded submodule of $M$, then $L \in \gr A$ also.
 
 \begin{prop}
\label{prop.fpres}
Let $A$ be a left graded-coherent ring with the property that every finite dimensional graded 
left $A$-module is 
finitely presented. Write $\cO$ for the image of $A$ in $\QGr A$. Then 
\begin{enumerate}
  \item 
  $\QGr A$ is a locally coherent category; 
  \item 
  the  full subcategory of $\QGr A$ consisting of the finitely presented objects is 
equivalent to $\qgr A$;
  \item 
$\cO$ is coherent: i.e.,  $\Hom_{\QGr A}(\cO,-)$ commutes with direct limits;
\item{}
 $\cO$ is finitely generated in the sense that whenever $\cO=\sum_i \cM_i$ for some directed family of subobjects $\cM_i$, $\cO=\cM_j$ for some $j$.
\end{enumerate}
\end{prop}
\begin{pf}
By definition, a Grothendieck category is locally coherent if its  full subcategory of finitely presented objects is abelian and every object is a direct limit of finitely presented objects.

Lazard \cite{Lazard} proved that every module over every ring  is a filtered direct limit of finitely presented modules (see, e.g., \cite[Example E.1.20]{Prest} or \cite{Osborne}). 
The proof for the ungraded case can be adapted to graded modules without difficulty. 
Thus every object in $\Gr A$ is a direct limit of objects in $\gr A$.

It follows that every object in $\QGr A$ is a direct limit of objects in $\qgr A$. 

Every object in $\Fdim A$ is a direct limit of objects in the Serre subcategory $\fdim A$ of $\gr A$ so 
condition (2) of \cite[Prop. A.4, p.112]{HK} is satisfied. Therefore  $\Fdim$ is 
a localizing subcategory of $\Gr A$ of finite type so    
\cite[Prop. A.5, p.113]{HK}  completes the proof. 
\end{pf}

\begin{remark}
The hypothesis in Proposition \ref{prop.fpres} that finite dimensional graded left $A$-modules are
finitely presented is satisfied if $A$ is a finitely generated $\NN$-graded algebra such that $\dim_k(A/A_{\ge 1}) < \infty$. For example, the hypothesis holds if $A$ is a graded quotient of a 
finitely generated free algebra whose generators have positive degree or, more generally, a 
a graded quotient of a path algebra of a finite quiver in which every arrow has positive degree. 
In particular, the graded algebras in this paper satisfy the hypotheses of Proposition \ref{prop.fpres}.
\end{remark}

\begin{lemma}
\label{lem.O.proj}
Let $\sA$ be a locally coherent category and ${\rm fp}\sA$ its full subcategory of finitely presented objects. Let $P \in {\rm fp}\sA$. If $P$ is projective as an object in ${\rm fp}\sA$, then it is 
is projective as an object in $\sA$.
\end{lemma}
\begin{pf}
Let $g:M \to N$ be an epimorphism in $\sA$ and suppose $f:P \to N$. 
We can write $M$ as a filtered direct limit of subobjects $M_i$ that belong to ${\rm fp}\sA$. Let $N_i=g(M_i)$. Then $N=\sum N_i$. It follows that $fP=\sum(fP \cap N_i)$. Since $P$ is finitely generated so is $fP$; therefore $fP= fP \cap N_j$ for some $j$. But that implies $fP \subset N_j$. Because $P$ is projective in ${\rm fp}\sA$, there is a morphism $h:P \to M$ with its image in $M_j$ such that $gh = f$. Hence $P$ is projective in $\sA$. 
\end{pf}

It is well-known that a free algebra is left and right coherent.  A general result, Lemma \ref{lem.coh}, then 
shows that $B$ is left and right coherent.

The next result is probably ``folklore''. 

\begin{lemma}
\label{lem.coh}
Let $A \subset B$ be compatibly $\ZZ$-graded $k$-algebras such that $B$ is a finitely 
presented left $A$-module. Suppose $A$ is left graded-coherent.  
Let $M \in \Gr B$. Then $M$ is finitely presented as a $B$-module if and only
if it is finitely presented as an $A$-module. 
In particular, $B$ is left graded-coherent.
\end{lemma}
\begin{pf}
($\Rightarrow$)
Let $0 \to K \to L \to M \to 0$ be an exact sequence of finitely generated graded left $B$-modules. Since $B$ is a finitely generated left $A$-module $K$ and $L$ are  finitely generated $A$-modules. Hence $M$ is a finitely presented $A$-module.

($\Leftarrow$)
Since $M$ is finitely generated as an $A$-module it is finitely generated as a $B$-module. Hence there 
is an integer $n$ and an exact sequence $0 \to K \to B^n \to M \to 0$ of graded left $B$-modules. By hypothesis, $B$, and therefore $B^n$, is a finitely presented left $A$-module. Since $A$ is left graded-coherent, 
$K$ is a finitely presented $A$-module. In particular, $K$ is finitely generated as a $B$-module. Hence $M$ is finitely presented as a $B$-module. 

Now we prove the final sentence of the lemma. 
Let $I$ be a finitely generated graded left ideal of $B$. Then $I$ is a finitely generated left $A$-module and, 
since $B$ is in $\gr A$, $I$ is also in $\gr A$. Therefore $I$ is finitely presented as a $B$-module.  
Hence $B$ is left graded-coherent.
\end{pf}

\begin{cor}
If $A$ and $B$ are as in Lemma \ref{lem.coh}, then the induction functor $B \otimes_A -:\Gr A \to \Gr B$
sends $\gr A$ to $\gr B$ and the restriction functor $\Gr B \to \Gr A$ sends $\gr B$ to $\gr A$.  
\end{cor}

\subsection{Graded coherent versus coherent}

This section, which plays no role in the other results in this paper, addresses a question raised by
the referee: what is known about the relation between graded coherence and coherence?

\begin{prop}
\label{prop.ref.qu}
Let $S=k[x,y]$ be the commutative polynomial ring in two variables and $R=k+xS$. Then $R$ is not coherent but is graded coherent with respect to the grading $\deg(x)=0$ and $\deg(y)=1$.\footnote{Tom Marley informed the author that he already knew that $R$ is not coherent.}
\end{prop}
\begin{pf}
To prove $R$ is graded coherent we use the following criterion of Piontkovski \cite[Prop. 3.2]{Pion1}:
{\it  An $\NN$-graded ring $A$ is right graded coherent if it has a right noetherian quotient ring 
$A/J$ such that $J$  is a free graded left $A$-module.}

Let $k[t]$ be the polynomial ring generated by a degree-one indeterminate, $t$.
Let $\theta:R \to k[y]$ be the restriction to $R$ of the homomorphism $\widetilde{\theta}:S \to k[y]$ defined by $\widetilde{\theta}(x)=1$ and $\widetilde{\theta}(y)=y$. The homogeneous components of $S$ are $S_0=k[x]$ and $S_n=k[x]xy^n$ for $n \ge 1$. Since $\theta$ preserves degree $\ker(\theta)$ is the direct sum of its homogeneous components: these are $\ker(\theta)_0=(x-1)k[x]$ and $\ker(\theta)_0=(x-1)k[x]xy^n$ for $n \ge 1$. Therefore $\ker(\theta)=(x-1)R$ which is a free graded $R$-module. 
 The image of $\theta$ is $k[t]$ so Piontkovski's criterion tells us that $R$ is graded coherent. 
 
 We will show that the ideal $I=xR+xyR$ is not finitely presented. The proof uses the 
$\ZZ \times \ZZ$-grading on $R$ that is inherited from that on $S$ given by
$$
\deg(x)=(1,0)
\qquad \hbox{and} \qquad 
\deg(y)=(0,1).
$$
Thus $S_{pq}=kx^py^q$.

The ideal $I$ is homogeneous, generated by elements of degree $(1,0)$ and $(1,1)$.
As usual, $R(m,n)$ denotes the graded $R$-module that is $R$ as an $R$-module and homogeneous 
components $R(m,n)_{p,q}=R_{p+m,q+n}$. 

Define the surjective degree-preserving homomorphism 
$$
\theta:F:=R(-1,0) \oplus R(-1,-1) \to I, \qquad \theta(a,b)=ax-bxy.
$$
The kernel, $K=\ker(\theta)$, is a graded submodule of $R(-1,0) \oplus R(-1,-1)$. 

The elements $(x^{p-1}y^q,0)$ and $(0,x^{p-1}y^{q-1})$ are a basis for $F_{pq}$ if $p \ge 2$ and 
$q\ge 1$. For example, $(xy,0)$ and $(0,x)$ are a basis for $F_{2,1}$.

We note that $(xy^j,xy^{j-1}) \in K_{2,j}$ for all $j \ge 1$ and $K_{p,q}=0$ for all $p \le 1$ and $q \le 0$.
Since $\fm =R_{(\ge 1, \ge 0)}$, it follows that the images in $K/\fm K$ of the elements  
$$
(xy^j,xy^{j-1}), \quad j \ge 1,
$$
are linearly independent. Thus, $\dim_k(K/\fm K) = \infty$. It follows that $K$ is not a finitely generated $R$-module. 
\end{pf}

The homogeneous components of $R$ for the grading $\deg(x)=0$ and $\deg(y)=1$ are not finite dimensional so, in some sense, the example in Proposition \ref{prop.ref.qu} is unsatisfying. It would be 
better to have an example of a non-coherent ring that is graded coherent with finite dimensional homogeneous components (if such an example exists!).

\subsection{Properties of $\QGr B$}  
Let $\sA$ be an abelian category, $\cO$ an object in $\sA$, and $(1)$ an automorphism of $\sA$ denoted by $\cF \mapsto \cF(1)$ on objects. For all $n \in \ZZ$, $\cF(n)$ has the obvious meaning.
 
 \begin{lemma}
 \label{lem.O.genor}
With the above notation, suppose 
$$
\cO \cong \cO(-1) \oplus  \bigoplus_{i=1}^p \cO(-r_i)
$$
for some integers $r_i >0$, and every $\cF \in \sA$ is a quotient of
$$
\bigoplus_{i=1}^q \cO(m_i) 
$$
for some integers $m_i$. Then for every $\cF \in \sA$ there is an epimorphism 
$$
\cO^n \twoheadrightarrow \cF
$$
for some integer $n$ that depends on $\cF$.
\end{lemma}
\begin{pf}
Let 
$$
\cG \cong \bigoplus_{i=1}^p \cO(-r_i).
$$

In this proof the notation $\cM  \twoheadrightarrow  \cN$ means ``there is an epimorphism from $\cM$ to $\cN$''. By hypothesis, $\cO  \twoheadrightarrow \cO(-1)$. An induction argument and repeated application of $(-1)$ shows that $\cO  \twoheadrightarrow \cO(-m)$ for all $m \ge 0$. Hence 
$\cO^p  \twoheadrightarrow \cG(1)$.
By hypothesis, $\cO(1) \cong \cO \oplus \cG(1)$ so $\cO^{p+1} \twoheadrightarrow \cO(1)$. An induction argument and repeated application of $(1)$ shows that for all $m\ge 0$, 
$\cO^q \twoheadrightarrow \cO(m)$ for some $q$ depending on $m$. Thus for every integer $m$, 
$\cO(m)$ is a quotient of a finite direct sum of copies of $\cO$. Hence every $\cF \in \sA$ is a quotient of a finite direct sum of copies of $\cO$.
\end{pf}

 We refer the reader to section \ref{intro.main} for definitions and notation concerning $\QGr B$ and $\qgr B$. In particular, $\b^*:\Gr B \to \QGr B$ denotes the quotient functor and $\cO=\b^*B$. 

\begin{prop}
\label{prop.O.genor}
Let $B$ be a connected graded $k$-algebra that is left graded coherent. Suppose $B$ contains a
left ideal of finite codimension that is isomorphic to $B(-1) \oplus G$ for some $G \in \gr B$. Then
\begin{enumerate}
  \item 
  Every object in $\qgr B$ is a quotient of $\cO^n$ for $n \gg 0$.
  \item 
  Every object in $\QGr B$ is a quotient of a direct sum of copies of $\cO$.  
  \item 
  $\cO$ is a generator in $\QGr B$.
\end{enumerate}
\end{prop}
\begin{pf}
Let $\cG=\b^* G$. Since $G$ is finitely generated, $\cG$ is a quotient of a finite direct sum of various $\cO(-i)$s. The hypothesis implies that $\cO \cong \cO(-1) \oplus \cG$. 
Lemma \ref{lem.O.genor} now implies that every $\cO(j)$ is a quotient of $\cO^n$ for a suitable $n$. 

(1) and (2).
Let $M \in \Gr B$ and let $\cM=\b^*M \in \QGr B$. Since $M$ is a quotient of a direct sum of various 
$B(m_i)$s,  $\cM$  is a quotient of a direct sum of various $\cO(m_i)$s. But each $\cO(m_i)$ is a quotient of a direct sum of copies of $\cO$. This completes the proof of (2). If $M$ is in $\gr B$, it   is a quotient of a direct sum of finitely many $B(m_i)$s so $\cM$ is a quotient of a finite direct sum of copies of $\cO$.
 
 (3)
Let  $f:\cM \to \cN$ be a non-zero morphism in $\QGr B$. Since $\cM$ is a quotient of a direct sum of copies of $\cO$ the restriction to one of the summands yields a non-zero composition $\cO \to \cM \stackrel{f}{\longrightarrow} \cN$. 
\end{pf}

Let $B$ be a connected graded $k$-algebra. 
Let 
$$
\sF:=\{\hbox{non-zero graded left ideals that are free and have finite codimension in $B$}\}.
$$

\begin{prop}
\label{prop.O.proj}
Let $B$ be a connected graded $k$-algebra. 
Suppose every
graded left ideal of finite codimension in $B$ contains some $I$ belonging to $\sF$.
Then $\cO$ is a projective object in $\QGr B$.
\end{prop}
\begin{pf}
By \cite[Cor.1, p.368]{Gab}, every short exact sequence in $\QGr B$  is of the form
$$
0 \longrightarrow  \b^* L \stackrel{\b^*f}{\longrightarrow}  \b^* M \stackrel{\b^*g}{\longrightarrow}  \b^* N \longrightarrow 0
$$
where $0 \longrightarrow  L \stackrel{f}{\longrightarrow}  M \stackrel{g}{\longrightarrow}   N \longrightarrow 0$ 
is an exact sequence in $\Gr B$. Let $\eta:\cO \to \b^*N$. By definition,
$$
\eta \in \liminj \Hom_{\QGr B}(B_{\ge i}, N/N') 
$$
where the direct limit is taken over all $i$ and all submodules $N' \subset N$ 
that are sums of finite dimensional submodules. The hypothesis implies that 
$$
\liminj \Hom_{\QGr B}(B_{\ge i}, N/N') = \liminj \Hom_{\QGr B}(I, N/N') 
$$ 
where the right-hand direct limit is taken over all $I \in \sF$ and all submodules $N' \subset N$ 
that are sums of finite dimensional submodules. 

Hence $\eta=\b^* h$ for some $h:I \to N/N'$ where $I\in \sF$. But $I$  is free so $h$ 
factors through $g$. Hence there is a morphism $\gamma:\cO \to \b^* M$ such that $\eta=(\b^*g) \circ \gamma$.  
\end{pf}
 
Let $B$ be any $\NN$-graded $k$-algebra.  
 If $f \in \End_{\Gr B} (B_{\ge n})$, then $f(B_{\ge n+1}) \subset B_{\ge n+1}$ so there is a directed system
\begin{equation}
\label{eq.dsys.End}
  \UseComputerModernTips
\xymatrix{
\cdots \ar[r] & \End_{\Gr B}(B_{\ge n})   \ar[r]^{\theta_n} & \End_{\Gr B}(B_{\ge n+1}) \ar[r] & \cdots
}
\end{equation}
of $k$-algebras in which $\theta_n(f) = f|_{B_{\ge n+1}}$.

\begin{thm}
\label{thm.progenor}
Suppose $B$ satisfies the hypotheses in Propositions  \ref{prop.O.genor} and \ref{prop.O.proj}. 
\begin{enumerate}
  \item 
  $\cO$ is a progenerator in $\QGr B$.
  \item 
The functor $\Hom(\cO,-)$ is an equivalence from the category $\QGr B$  to the category of right modules over the endomorphism ring $\End_{\QGr B} \cO$.
\item{}
If every $B_{\ge n}$ contains a non-zero free right $B$-module, then 
$\End_{\QGr B} \cO \cong \liminj \End_{\Gr B}(B_{\ge n})$, the direct limit of (\ref{eq.dsys.End}). 
\end{enumerate}
\end{thm}
\begin{pf}
Part (1) is  the content of  Propositions  \ref{prop.O.genor} and \ref{prop.O.proj} and part (2)
is a standard consequence of (1).

(3)
 By the definition of  morphisms in a quotient category,
\begin{equation}
\label{.d.sys1}
\End_{\QGr B}\cO =  \liminj \Hom_{\Gr B}(B',B/B'')
\end{equation}
where $B'$ runs over all graded left ideals in $B$ such that $\dim_k(B/B') < \infty$ and 
$B''$ runs over all graded left ideals in $B$  
that are sums of finite dimensional left ideals. 

Let $B''$ be a finite dimensional graded left ideal in $B$.
Then $B_{\ge n}B''=0$ for some $n$. The hypothesis in (3) implies that  
$B''=0$. Since the system of left ideals $B'$  such that
 $\dim_k(B/B') < \infty$  is cofinal with the system of left ideals
$B_{\ge i}$,  
$$
\End_{\QGr B} \cO =\liminj_{i} \Hom_{\Gr B}(B_{\ge i},B).
$$
But the map $ \Hom_{\Gr B}(B_{\ge i},B_{\ge i}) \to \Hom_{\Gr B}(B_{\ge i},B)$ induced by 
the inclusion $B_{\ge i} \to B$ is an isomorphism so
$$
\End_{\QGr B} \cO = \liminj_{i} \Hom_{\Gr B}(B_{\ge i},B_{\ge i}).
$$
 This completes the proof.
\end{pf}

\section{The algebra $k \langle x,y\rangle/(y^{r+1})$}

In this section $r$ is a  positive integer and
$$
B:=\frac{k\langle x,y\rangle}{(y^{r+1})}.
$$

\subsection{The subalgebra $A$}
\label{ssect.A+B}
Let $A$ be the subalgebra of $B$ generated by $u_1,\ldots,u_{r+1}$ where
$$
u_j:= \sum_{i=1}^j y^{i-1}xy^{j-i}, \qquad 1 \le j \le r+1.
$$
Thus, $u_1=x$, $u_2=xy+yx$, $u_3=xy^2+yxy+y^2x$, and so on. 

The {\sf Hilbert series} of $B$ is $H_B(t)=\sum_{i\ge 0} b_it^i$ where 
$$
b_i=\dim_k B_i.
$$

\begin{prop}
\label{prop.B}
Let $A$ and $B$ be as above. 
\begin{enumerate}
  \item 
  $A$ is the free $k$-algebra on $u_1,\ldots,u_{r+1}$.
  \item{}
  $B$ is free as a left, and as a right, $A$-module with basis $\{1,y,\ldots,y^r\}$.
  \item{}
  The Hilbert series for $B$ is  
  $$
  H_B(t) =\frac{1+t+\cdots +t^r}{1-t-\cdots -t^{r+1}}.
  $$
 \item 
  $H_A(t)=(1-t -\cdots -t^{r+1})^{-1}$.
  \item{}
  $b_i=2^i$ for $0\le i \le r$, and 
  $$
  b_{n+1}=b_n+b_{n-1}+\cdots + b_{n-r}
  $$
  for $n \ge r$. 
\end{enumerate}
\end{prop}
\begin{pf}
(3) and (6).
A basis for $B_n$, the degree $n$ component of $B$, is provided by
$$
W_n:=\{\hbox{words of length $n$ in $x$ and $y$ not containing $y^{r+1}$ as a subword}\}.
$$
We introduce the following notation
\begin{align*}
W_{in}&=\{w \in W_n \; | \; w=w'xy^i \; \hbox{for some } \, w' \in W_{n-i-1}\}, \quad 0 \le i \le r,
\\
b_{in}&=|W_{in}|. 
\end{align*}
Since $W_n$ is the disjoint union of the $W_{in}$s, $b_n=b_{0n}+\cdots + b_{rn}$.

Since 
$
W_{i+1,n+1}=\{wy \; | w  \in W_{in}\} = W_{in}y,
$
we have $b_{i+1,n+1}=b_{in}$. Since
$
W_{0n}=\{wx \; | w \in W_{n-1}\} = W_{n-1}x,
$
we have $b_{0n}=b_{n-1}$. Therefore
$$
b_{n+1}=b_n+b_{n-1}+\cdots + b_{n-r}
$$
for all $n \ge r$.  Also $b_i=2^i$ for $0 \le i \le r$. (When $r=1$ we have $b_0=1$, $b_1=2$, and 
$b_{n+1}=b_n+b_{n-1}$, so $(b_n)_{n \ge 0}$ is the Fibonacci sequence beginning $1,2,3,5,\ldots$. When $r=2$, the sequence $(b_n)_{n \ge 0}$ is called the Tribonacci sequence. The limit of the ratio of successive terms $b_{n}/b_{n+1}$ of the Tribonacci sequence is a zero of the polynomial $1-x-x^2-x^3$. The Tribonacci sequence is related to the Rauzy fractal \cite{R}.)
Multiplying both sides by $t^{n+1}$ and summing over all $n \ge r$ we obtain
$$
\sum_{n=r+1}^\infty b_nt^n = \sum_{j=0}^{r} t^{j+1} \sum_{n=r-j}^\infty b_nt^n.
$$
In other words,  
$$
H_B(t)-\sum_{i=0}^r(2t)^i =  \sum_{j=0}^{r} t^{j+1} \Bigg(H_B(t) -  \sum_{i=0}^{r-j-1} (2t)^i\Bigg)
$$
so
$$
(1-t-t^2-\cdots-t^{r+1})H_B(t) = \sum_{i=0}^r(2t)^i - \sum_{j=0}^{r-1} t^{j+1} \sum_{i=0}^{r-j-1} (2t)^i
$$
By induction on $r$, the case $r=1$ being clear, one can show that the right-hand side of this equation is equal to $1+t+\cdots + t^r$. Hence the formula for $H_B(t)$ in (3) is correct.

We now turn to the proof of parts (1), (2), and (4), of the proposition.

The first step is to show that $M:=A+Ay+\cdots +Ay^r$ is equal to $B$  by showing it is a right 
ideal of $B$ (that contains 1). Clearly, $My \subset M$ so it remains to show $Mx \subset M$. The key calculation is
$$
u_j-u_{j-1}y=\sum_{i=1}^j y^{i-1}xy^{j-i} - \sum_{i=1}^{j-1} y^{i-1}xy^{j-1-i}y = y^{j-1}x
$$
for $2 \le j \le r+1$. This implies that $Ay^{j-1}x \subset Au_j+Au_{j-1}y$; but $u_j$ and $u_{j-1}$ belong
to $A$ so 
$$
Ay^{j-1}x \subset A+Ay
$$
if $2 \le j \le r+1$. Since $x \in A$, $Ax \subset A$ so we conclude that $Mx \subset M$. Hence
$$
A+Ay+\cdots +Ay^r=B
$$
as claimed. A similar argument shows that $A+yA+\cdots +y^rA=B$.

We now consider the free algebra $F=k\langle v_1,\dots, v_{r+1} \rangle$ with grading given by 
$\deg v_i=i$. Since $F$ is the free coproduct $k[v_1]*k[v_2]*\cdots *k[v_{r+1}]$ of polynomial rings, the Hilbert series formula
$$
H_{R*S}(t)^{-1} = H_R(t)^{-1} +H_S(t)^{-1}-1
$$
gives
$$
H_F(t)=(1-t-\cdots -t^{r+1})^{-1}.
$$
The homomorphism $\psi:F \to A$ given by $\psi(v_i)=u_i$ allows us to view $B$ as a graded $F$-module. We have just shown that $F+Fy+\cdots+Fy^r =B$ so there is a surjective homomorphism
of graded $F$-modules 
$$
\Psi:F \oplus F(-1) \oplus \cdots \oplus F(-r) \twoheadrightarrow B.
$$
However, the Hilbert series of this free $F$-module is $(1+t+\cdots +t^r)H_F(t)$ which is equal to 
$H_B(t)$ so   $\Psi$ must be an isomorphism. 

It follows at once that $\psi:F \to A$ is an isomorphism, thus proving (1) and that $B$ is a free $A$-module on the right and on the left with basis $\{1,y,\ldots,y^r\}$, thus proving (2). The calculation of $H_F(t)$ and the isomorphism $F \cong A$ proves (4).
\end{pf}

\begin{prop}
\label{prop.B.coh}
The algebra $B=k\langle x,y\rangle/(y^{r+1})$ is left and right coherent.
\end{prop}
\begin{pf}
This follows from Lemma \ref{lem.coh} and Proposition \ref{prop.B}.
\end{pf}

\begin{lemma}
\label{lem.free.left.ideals}
Let  $B=k\langle x,y\rangle/(y^{r+1})$ and let $A$ be its subalgebra generated by $u_1,\ldots,u_{r+1}$.
Then
\begin{enumerate}
\item 
  $BA_{\ge 1} \cong B(-1) \oplus \cdots \oplus B(-r-1)$;
  \item 
  $BA_{\ge n}$ is a free left $B$-module for all $n \ge 0$;
  \item{}
  $\dim_k (B/BA_{\ge n}) <\infty$ for all $n \ge 0$;
  \item 
  $BA_{\ge 1}=BxB$ and $B/BA_{\ge 1} \cong k[y]/(y^{r+1})$.
\end{enumerate}
\end{lemma}
\begin{pf}
(1)
Since $A$ is a free algebra every left ideal in $A$ is a free $A$-module. In particular, 
$$
A_{\ge 1}= Au_1+\cdots+Au_{r+1} \cong A(-1) \oplus \cdots \oplus A(-r-1).
$$
Since $B$ is a free right $A$-module, applying the functor $B \otimes_A -$ to the exact sequence
$0 \to A_{\ge 1} \to A \to A/A_{\ge 1} \to 0$ gives an exact sequence
\begin{equation}
\label{imp.ses}
0 \to \bigoplus_{i=1}^{r+1} B(-i) \to B \to B/BA_{\ge 1} \to 0.
\end{equation}
This proves (1).

(2)
This is proved in the same way as (1). Because $A$ is a free algebra, for every $n \ge 0$,
there is a graded vector space $V$, depending on $n$, and an exact sequence
$0 \to A \otimes_k V \to A \to A/A_{\ge n} \to 0$.  After applying  $B \otimes_A -$ to this sequence
we see that $BA_{\ge n} \cong B \otimes_k V$.

(3)
Certainly $A/A_{\ge n}$ has finite dimension. Since $B$ is a finitely generated right $A$-module it follows that $B \otimes_A(A/A_{\ge n}$ has finite dimension. But $B \otimes_A(A/A_{\ge n} \cong B/BA_{\ge n}$. 

(4)
Set $u_0=0$. Then $xy^{j-1}=u_j - yu_{j-1}$ for all $1 \le j \le r+1$. Hence
\begin{equation}
\label{eq.Buj}
\sum_{j=1}^{r+1} Bxy^{j-1} \subset \sum_{j=1}^{r+1} Bu_j.
\end{equation}
On the other hand, $Bu_1=Bx$ and, since $u_j=xy^{j-1} +yu_{j-1}$, $Bu_j \subset Bxy^{j-1}+Bu_{j-1}$
for all $2 \le j \le r+1$. An induction argument shows that 
$$
Bu_1+\cdots +Bu_j \subset Bx+Bxy+\cdots Bxy^{j-1}
$$
for all $1 \le j \le r+1$. The case $j=r+1$, together with (\ref{eq.Buj})  shows that  
\begin{equation}
\label{eq.Buj.2}
\sum_{j=1}^{r+1} Bu_j = \sum_{j=1}^{r+1} Bxy^{j-1}.
\end{equation}
The right-hand sum is closed under right multiplication by $y$; it is also closed under right-multiplication by $x$ because the first term in the sum is $Bx$; it is therefore a two-sided ideal of $B$. The left-hand side of (\ref{eq.Buj.2})  is $BA_{\ge 1}$ because $A_{\ge 1}=Au_1+\cdots +Au_{r+1}$. Hence  $BA_{\ge 1}$ is
a two-sided ideal of $B$. 

Since $x \in A$,  $BxB \subset BA_{\ge 1}$. However, the right-hand side of 
(\ref{eq.Buj.2})
is contained in $BxB$ so $BA_{\ge 1} = BxB$. It is obvious that $B/(x) \cong k[y]/(y^{r+1})$. 
\end{pf}

\subsection{Remark}
It follows from parts (1) and (3) of Lemma \ref{lem.free.left.ideals} that
$$
\cO \cong \cO(-1) \oplus \cdots \oplus \cO(-r-1)
$$
in $\qgr B$. This fact, which has no parallel in projective algebraic geometry or noetherian ring theory,
tells us $\Projnc B$ is a very unusual geometric object when compared to those that arise in 
classical algebraic geometry.  

\subsection{Remark} 
Michel Van den Bergh pointed out that when $r=1$ there are  isomorphisms
$$
\cO \cong \cO(-1) \oplus \cO(-2) \cong \cdots \cong \cO(-n)^{f_n}\oplus \cO(-n-1)^{f_{n-1}}
$$
for all $n \ge 1$. 

\section{Proof of Theorem \ref{thm.main}}
   
\subsection{}
In this section $B=k\langle x,y\rangle/(y^{r+1})$.

\begin{prop}
Let $B$ be as above.
\begin{enumerate}
  \item 
  $\cO$ is a progenerator in $\QGr B$.
  \item 
The functor $\Hom(\cO,-)$ is an equivalence from $\QGr B$  to the category of right modules 
over the endomorphism ring $\End_{\QGr B} \cO$.
\item{} 
$\End_{\QGr B} \cO \cong \liminj \End_{\Gr B}(B_{\ge n})$, the direct limit of (\ref{eq.dsys.End}). 
\end{enumerate}
\end{prop}
\begin{pf}
By Lemma \ref{lem.free.left.ideals},  $BA_{\ge n}$ has finite codimension in $B$ and is isomorphic to
$B(-1) \oplus \cdots \oplus B(-r-1)$. Hence $B$ satisfies the hypothesis of Proposition \ref{prop.O.genor}.
 
 If $I$ is a graded left ideal of finite codimension in $B$, then $I$ contains $B_{\ge i}$ for some $i$. 
By Lemma \ref{lem.free.left.ideals},  
$BA_{\ge i}$ is a free $B$-module contained in $B_{\ge i}$ and has finite codimension in $B$  
so $B$ satisfies  the hypotheses of Proposition \ref{prop.O.proj}. Thus, (1) and (2) follow from
Theorem \ref{thm.progenor}.

Since $B$ is anti-isomorphic to itself, the fact that $B_{\ge n}$ contains a non-zero free left $B$-module 
implies it contains a non-zero free right $B$-module. Hence $B$ satisfies the hypothesis in part (3) of Theorem \ref{thm.progenor} and therefore its conclusion which is part (3) of the present proposition.
\end{pf}

Our next task is to give a precise description of the directed system (\ref{eq.dsys.End}) and its direct limit. Lemma \ref{lem.main} shows that each term, $\End_{\Gr B}(B_{\ge n})$, in the directed system is a block lower triangular matrix algebra that is 
Morita equivalent to the ring of lower triangular $(r+1)\times(r+1)$ matrices over $k$. Lemma \ref{lem.phin.Sn} shows that the diagonal, or semisimple,  part of each $\End_{\Gr B}(B_{\ge n})$
forms a directed system that has the same direct limit as (\ref{eq.dsys.End}). Thus $\End_{\QGr B}\cO$ is 
the direct limit of a directed system in which each term is a product of $r+1$ matrix algebras.

Define  $B_n^i:=y^ixB_{n-i-1}$ for $0 \le i \le r$. Then
$$
B_n=B^0_n \oplus B^1_n \oplus \cdots \oplus B^r_n.
$$

\begin{lem}
\label{lem.main}
The restriction map 
$$
\Phi_n:\End_{\Gr B} B_{\ge n} \to \End_k B_{n}, 
\qquad 
\Phi_n(f):=f|_{B_n},
$$
is an isomorphism from $\End_{\Gr B} B_{\ge n}$ onto the lower triangular subalgebra
$$
T_n  : =\begin{pmatrix}
 \Hom_k(B^0_n,B^0_n)     &  0  & 0 & \cdots & 0 \\
 \Hom_k(B^0_n,B^1_n)     &   \Hom_k(B^1_n,B^1_n)   & 0 & \cdots & 0 \\
 \vdots & \vdots  & \vdots &  & \vdots \\
  \Hom_k(B^0_n,B^r_n)     &   \Hom_k(B^1_n,B^r_n)   &   \Hom_k(B^2_n,B^r_n) & \cdots &   \Hom_k(B^r_n,B^r_n) 
\end{pmatrix}
$$
of $\End_k B_n$.
\end{lem}
\begin{pf}
Because $B$ is generated as an algebra by $B_1$, $B_{\ge n}$ is generated as a left ideal by $B_n$.
Hence $\Phi_n$ is injective. 
 
Let $f \in \End_{\Gr B} B_{\ge n}$.  
Since $y^{r-i+1}B^i_n=0$, $y^{r-i+1}f(B^i_n)=0$.  Therefore
$$
f(B^i_n) \subset B^i_n \oplus B^{i+1}_n \oplus \cdots \oplus B^r_n.
$$
Hence $\Phi_n(f)$ belongs to $T_n$. 

Fix $i$ and suppose $g \in T_n$ is such that $g(B^j_n)=0$ for all $j \ne i$. 
It is clear that
$$
B_{\ge n}=BB^0_n \oplus BB^1_n \oplus \cdots \oplus BB^r_n.
$$
Let $f:B_{\ge n} \to B_{\ge n}$ be the unique linear map that is zero on $BB^j_n$ for all $j \ne i$ and
satisfies $f(ab)=ag(b)$ if $a \in B$ and $b \in B^i_n$.
The fact that $g(B^i_n) \subset B^i_n \oplus B^{i+1}_n \oplus \cdots \oplus B^r_n$ ensures that $f$ 
is a well-defined left $B$-module homomorphism, and the condition $f(ab)=ag(b)$ ensures that 
$\Phi_n(f)=g$.

Now let $g$ be an arbitrary element of $T_n$. Then $g=g_0+g_1+\cdots +g_r$ where $g_i$ is defined
to be $g$ on $B^i_n$ and zero on $B^j_n$ for all $j \ne i$. Let $f_i :B_{\ge n} \to B_{\ge n}$ be the unique linear map that vanishes on $BB^j_n$ for $j \ne i$ and satisfies 
$f_i(ab)=ag_i(b)$ if $a \in B$ and $b \in B^i_n$. Now $f_0+f_1+\cdots + f_r$ is in $\Hom_{\Gr B}(B_{\ge n}, B_{\ge n})$ and 
$\Phi_n(f_0+f_1+\cdots + f_r)=g$. 
This completes the proof that the image of $\Phi_n$ is $T_n$.
\end{pf}

Let 
$$
\theta_n:\End_{\Gr B}(B_{\ge n}) \to \End_{\Gr B}(B_{\ge n+1})
$$
be the restriction map. Let $f \in \Hom_{\Gr B}(B_{\ge n},B_{\ge n})$. The bottom arrow in the diagram
\begin{equation}
\label{comm.diag.Tn}
  \UseComputerModernTips
\xymatrix{
\End_{\Gr B}(B_{\ge n}) \ar[d]_{\Phi_n} \ar[rr]^{\theta_n} && \End_{\Gr B}(B_{\ge n+1}) \ar[d]^{\Phi_{n+1}}
\\
T_n \ar[rr]_{\psi_n:=\Phi_{n+1} \circ \theta_n \circ \Phi_n^{-1}} && T_{n+1}.
}
\end{equation}
is  $\psi_n\big(\Phi_n(f)\big)=\psi_n\big(f|_{B_n}\big)=f|_{B_{n+1}}$.   
  
\subsection{The algebras $S_n$}  
Let $S_n$ be the diagonal subalgebra of $T_n$, i.e., 
 $$
S_n : =\begin{pmatrix}
 \Hom_k(B^0_n,B^0_n)     &  0  & 0 &  & 0 \\
0     &   \Hom_k(B^1_n,B^1_n)   & 0 &  & 0 \\
 \vdots & \vdots  &  \ddots && \vdots \\
  0    &  0   &   \cdots &  &   0\\
  0    &  0   &   \cdots &  &   \Hom_k(B^r_n,B^r_n) 
\end{pmatrix}
$$ 
Since $\dim_k B^i_n=\dim_k y^ixB_{n-i-1}=\dim_k B_{n-i-1} = b_{n-i-1}$,  
$S_n$ is a product of $r+1$ matrix algebras of sizes $b_{n-1},b_{n-2},\ldots,b_{n-r-1}$. 

\begin{lemma}
\label{lem.phin.Sn}
Let $(b_0,\ldots,b_r) \in B_{n+1}^0 \times \cdots \times B_{n+1}^r$ and 
$$
g=(g_0,\ldots,g_r) \in \big(\End_k B_n^0\big) \times \cdots \times  \big(\End_k B_n^r\big) \, = \, S_n.
$$
\begin{enumerate}
  \item 
  There are unique elements $c_j,d_j \in B_n^j$ such that
  $$
  (b_0,\ldots,b_r) = \Big(x\sum_{j=0}^r c_j,yd_0,yd_1,\ldots,yd_{r-1}\Big).
  $$
  \item 
  There is a well-defined map $\phi_n:S_n \to S_{n+1}$ defined by
  \begin{align*}
  \phi_n(g_0,\ldots,g_r)(b_0,\ldots,b_r) := &
  \\
    \Bigg(x\sum_{j=0}^r g_j(c_j),& \, yg_0(d_0),yg_1(d_1),\ldots,yg_{r-1}(d_{r-1})\Bigg)
  \end{align*}
  where the $c_j$s and $d_j$s are as in (1).
  \item 
  If $\a_n:S_n \to T_n$ is the inclusion, then the diagram
  $$
  \UseComputerModernTips
\xymatrix{
S_n \ar@<.6ex>[d]_{\a_{n}} \ar[rr]^{\phi_n} && S_{n+1}  \ar@<.6ex>[d]^{\a_{n+1}} 
\\
 T_n  \ar[rr]_{\psi_n}  && T_{n+1}  
}
$$
commutes.
\item{}
There is a homomorphism of directed systems $\a_{\ldot}:(S_n,\phi_n) \longrightarrow (T_n,\psi_n)$.
\end{enumerate}
\end{lemma}
\begin{pf}
(1)
Because $b_0 \in B_{n+1}^0=xB_n=x(B_n^0 \oplus B_n^1 \oplus \cdots \oplus B_n^r)$ and $x$ is not a zero divisor there are unique elements $c_j \in B_n^j$ such that $b_0=x(c_0+\cdots +c_r)$.

Suppose $1 \le i \le r$. Then $b_i \in B_{n+1}^i=y^ixB_{n-i}=yB_n^{i-1}$ and because left multiplication 
by $y$ is a bijective map $B_n^{i-1} \to B_{n+1}^i$ there is a unique $d_{i-1} \in B_n^{i-1}$ 
such that $b_i=yd_{i-1}$. 

(2)
Since $c_j,d_j \in B_n^j$, $g_j(c_j)$ and $g_j(d_j)$ belong to $B_n^j$ too. Hence $xg_j(c_j) \in B_{n+1}^0$ and $yg_j(d_j) \in yB_n^j =B_{n+1}^{j+1}$. Hence $\phi_n(g_0,\ldots,g_n) \in S_{n+1}$
and $\phi_n:S_n \to S_{n+1}$.

(3)
Since $T_{n+1} \subset \End_k B_{n+1}$ and $B_{n+1}=B_{n+1}^0 \oplus \cdots \oplus B_{n+1}^r$
we must show that
$$
\psi_n\a_n(g_0,\ldots,g_r)(b_0+\cdots+b_r)=\a_{n+1}\phi_n(g_0,\ldots,g_r)(b_0+\cdots+b_r)
$$
for all $(g_0,\ldots,g_r) \in S_n$ and $b_i \in B_{n+1}^i$, $0 \le i \le r$.

Define $f_i \in \End_{\Gr B}(B_{\ge n})$ to be the unique linear map such that
\begin{enumerate}
  \item 
  $f_i$ is zero on $BB^j_n$ for all $j \ne i$ and
  \item 
$f_i(ab)=ag_i(b)$ if $a \in B$ and $b \in B^i_n$.
\end{enumerate} 
The fact that $g_i(B^i_n) \subset B^i_n$ ensures that $f_i$ is a well-defined left $B$-module homomorphism and (2) ensures that $\Phi_n(f_i)=g_i$. Define $f=f_0+\cdots +f_r$. Then 
$\Phi_n(f)=\diag(g_0,\cdots,g_r)$. Therefore
$$
\psi_n\a_n(g_0,\ldots,g_r) = \Phi_{n+1}\theta_n(f) = f\big\vert_{B_{n+1}}
$$
and 
\begin{align*}
\psi_n\a_n(g_0,\ldots,g_r)&(b_0+\cdots+b_r) = f(b_0+\cdots+b_r)
\\
&= f\Big(x\sum_{j=0}^r c_j + y(d_0+\cdots +  d_{r-1} )\Big)
\\
&=x\sum_{j=0}^r g_j(c_j) + yg_0(d_0) +  \cdots + yg_{r-1}(d_{r-1} ).
\end{align*}
On the other hand, 
\begin{align*}
\a_{n+1}\phi_n(g_0,\ldots,g_r)&(b_0+\cdots+b_r) =  \diag\big(\phi_n(g_0,\ldots,g_r)\big)(b_0+\cdots+b_r) 
\\
&= \sum_{i=0}^r \phi_n(g_0,\ldots,g_r)(b_i)
\\
&=x\sum_{j=0}^r g_j(c_j) + yg_0(d_0) +  \cdots + yg_{r-1}(d_{r-1} )
\end{align*} 
where the last equality follows from part (2). This completes the proof that the diagram in (3) commutes
and (4) is an immediate consequence. 
\end{pf}

\begin{thm}
\label{thm.S=T}
The homomorphisms $\a_{n}:(S_n,\phi_n) \longrightarrow (T_n,\psi_n)$ in Lemma
\ref{lem.phin.Sn} induce an isomorphism
$$
\a: \liminj_n (S_n,\phi_n) \stackrel{\sim}{\longrightarrow}  \liminj_n  (T_n,\psi_n).
$$
\end{thm}
\begin{pf}
Write $S$ and $T$ for the direct limits of the two systems. Lemma \ref{lem.phin.Sn} showed that the
diagonal inclusions $\a_n :S_n \to T_n$ induce an algebra homomorphism $\a:S \to T$. Every $\a_n$ is 
injective, as are the $\phi_n$s and $\psi_n$s, so $\a$ is injective. It remains to show that $\a$ is surjective.

Let $t \in T_n$. Then $t=\Phi_n(f)$ for some  $f \in \End_{\Gr B}(B_{\ge n})$ and
 $$
\psi_{n+r-1} \ldots \psi_n(t)=   \psi_{n+r-1} \ldots \psi_n\Phi_{n}(f) = \Phi_{n+r} \theta_{n+r-1} \ldots \theta_n(f) = \Phi_{n+r}(f')
$$
where $f'$ is the restriction of $f$ to $B_{\ge n+r}$.

Let $0 \le i \le r-1$. Then $n+r-i-1 \ge n$ so 
$$
f(B^i_{n+r})=f(y^ixB_{n+r-i-1}) \subset y^ixf(B_{n+r-i-1}) \subset y^ixB_{n+r-i-1} = B^i_{n+r}.
$$
Also, 
$$
f(B^r_{n+r})=f(y^rxB_{n-1}) = y^rf(xB_{n-1}) \subset y^rB_n = y^rxB_{n-1} =B^r_{n+r}.
$$
This shows that $\Phi_{n+r}(f')$ belongs to $S_{n+r}$ or, more precisely, $\psi_{n+r-1} \ldots \psi_n(t)$
is in the image of $\a_{n+r}$. It follows that $\a$ is surjective and therefore an isomorphism. 
\end{pf}

\subsection{Definition of $S$}
\label{sect.S}
From now on we will write 
$$
S:=  \liminj_n (S_n,\phi_n)
$$
where the direct limit is taken in the category of $k$-algebras. 
Recall that $\Mod S$ is the category of right $S$-modules.

 \begin{thm}
There is an equivalence of categories 
$$
\QGr B \equiv \Mod S.
$$
Under the equivalence, $\cO$ corresponds to $S_S$.
\end{thm}
\begin{pf}
Because (\ref{comm.diag.Tn}) commutes the directed system $\big(\End_{\Gr B}(B_{\ge n}),
\theta_n\big)$ is isomorphic to the directed system $(T_n,\psi_n)$.
By Theorem \ref{thm.progenor},  Lemma \ref{lem.main}, and Theorem
\ref{thm.S=T},
$$
\QGr B \equiv \Mod \End_{\QGr B}\cO   \equiv \Mod \liminj (T_n,\psi_n) = \Mod \liminj_n (S_n,\phi_n).
$$
The equivalence is given by $\Hom_{\QGr B}(\cO,-)$ so sends $\cO$ to $S_S$.
\end{pf}

\section{Properties of $S$ and consequences for $\QGr B$} 

We refer the reader to \cite{Eff} for information about Bratteli diagrams.
 
 We write $S_n=S_n^0 \times \cdots \times S_n^r$ where $S_n^i$ is the matrix algebra of size 
$b_{n-i-1}$.
 
\begin{prop}
\label{prop.Brat.Sn}
The Bratteli diagram for $S$ is repeated copies of 
$$
  \UseComputerModernTips
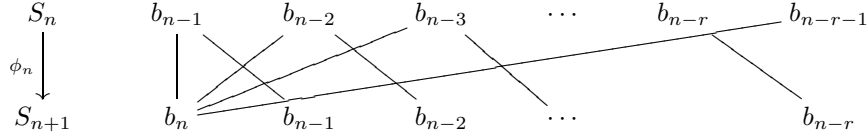
\xymatrix{
S_n \ar[d]_{\phi_n} & b_{n-1}  \ar@{-}[dr]  \ar@{-}[d] & b_{n-2} \ar@{-}[dl] \ar@{-}[dr]  & b_{n-3}   \ar@{-}[dr]   \ar@{-}[dll]  & \cdots & b_{n-r} \ar@{-}[dr]  & b_{n-r-1}  \ar@{-}[dlllll] 
\\
S_{n+1} & b_n & b_{n-1} &  b_{n-2}   & \cdots &   & b_{n-r} 
}
$$
where $b_i$ is $\dim_k B_i$ and is  given explicitly in Proposition \ref{prop.B}(5). 
\end{prop}
\begin{pf}
If we think of $S_n$ abstractly as a product of  $r+1$ matrix algebras of sizes 
$b_{n-1},b_{n-2},\ldots,b_{n-r-1}$,
and not concretely as endomorphisms of $B_n^0 \times \cdots \times B^r_n$, then the map 
$\phi_n:S_n \to S_{n+1}$ described in Lemma \ref{lem.phin.Sn}(2) is 
$$
\phi_n(g_0,\ldots,g_r)=\big(\diag(g_0,\ldots,g_r),g_0,g_1,\ldots,g_{r-1}\big).
$$
The result follows.
\end{pf}

\begin{prop}
$S$ is a simple von Neumann regular ring.
\end{prop}
\begin{pf}
 A ring $R$ is von Neumann regular if for each $x \in R$ the equation $xyx=x$ has a solution $y \in R$.
A matrix algebra  over a field is von Neumann regular so a product of matrix algebras is von Neumann regular. It follows from the definition that a filtered direct limit of von Neumann regular rings is von Neumann regular. Hence $S$ is von Neumann regular. 

The simplicity of $S$ can be read off from the shape of its Bratteli diagram. 
Let $x$ be a non-zero element of $S$.
Then $x \in S_n$ for some $n$ so can be written as $(x_0,\ldots,x_r)$ where $x_i \in S_n^i$.
Some $x_i$ is non-zero so, as can be seen from the Bratteli diagram, the image of $x$ in $S_{n+r}$
has a non-zero component in $S_{n+r}^j$ for all $j=0,\ldots,r$. Hence the ideal of $S_{n+r}$ generated by the image of $x$ is $S_{n+r}$ itself. Hence $SxS=S$.
\end{pf}

Other properties of $S$ can be read off from its description as a filtered direct limit of products of matrix algebras. For example, $S$ is unit regular \cite[Ch. 4]{vNR} hence directly finite \cite[Ch. 5]{vNR},
and it satisfies the comparability axiom \cite[Ch. 8]{vNR}.

\begin{prop}
The ring $S$ is coherent on both the left and right.
\end{prop}
\begin{pf}
Since matrix algebras are isomorphic to their opposite rings, $S$ is isomorphic to its opposite and 
it therefore suffices to prove the result for just one side.
We have already shown that $\QGr B$, and hence $\Mod S$, is a locally coherent category so $S$ is right coherent.

Alternatively, $S$ has global homological dimension one (see section \ref{ssect.curve}) and every ring of global dimension one is coherent; to see that let $I$ be a finitely generated left ideal and  $P$ a module such that $I \oplus P$ is isomorphic to a finitely generated free module, $F$ say; then $I \cong F/P$ and $P \cong F/I$ so $P$ is finitely generated and $I$ is therefore finitely presented.  
\end{pf}

Translating standard properties of von Neumann regular rings into statements about $\QGr B$ and 
$\qgr B$ shows these categories are very different from the categories $\Qcoh X$ and $\coh X$ that
appear in algebraic geometry.

\begin{prop}
\label{prop.hered}
Every subobject of $\cO$ is projective.
\end{prop}
\begin{pf}
Since $S$ is a countable union of matrix algebras it has countable dimension.
So every left ideal of $S$ is countably generated and therefore projective by \cite[Cor. 2.15]{vNR}. 
\end{pf}

\begin{prop}
Every object in $\qgr B$ is projective in $\QGr B$ and every short exact sequence in $\qgr B$ splits.
\end{prop}
\begin{pf}
Every module over a von Neumann regular ring is flat \cite[Cor. 1.13]{vNR}. Hence all finitely presented $S$-modules are projective. But the equivalence $\QGr B \equiv \Mod S$ restricts to an equivalence
$\qgr B \equiv \mod S$ so all objects in $\qgr B$ are projective in $\QGr B$. 
\end{pf}

 \begin{prop}
There are no noetherian objects in $\qgr B$ except 0. 
 \end{prop}
 \begin{pf}
It suffices to show that $\mod S$ does not contain a simple module.  Suppose to the contrary that $L$ is a finitely presented simple left $S$-module. There is a surjective homomorphism $S \to L$ and its kernel, $K$ say, is a finitely generated left ideal of $S$. Since $S$ is von Neumann regular $K=S(1-e)$ for some idempotent $e$ and $L \cong Se$. Since $L$ is simple $e$ is a primitive idempotent, but 
it follows from the Bratteli diagram that a primitive idempotent in $S_n$ is no longer 
primitive in $S_{n+2}$. Hence $S$ has no primitive idempotents.
 \end{pf}

\section{Multi-matrix algebras associated to partitioned sets}
\label{sect.inv.lim.dlim}

A finite product of finite dimensional matrix algebras over a common field is called a {\sf multi-matrix algebra}. 

In \cite[Sect. II.3]{NCG} Connes associates to the space of Penrose tilings a C${}^*$-algebra 
constructed as a direct limit of multi-matrix algebras. 
In this section we generalize his construction (the $r=1$ case is the algebra Connes constructs)
and show the ring $S$ in section \ref{sect.S} is isomorphic to the ring obtained by this generalization.

\subsection{An algebra homomorphism associated to a map between partitioned sets}

Let $X$ be a finite set  endowed with a partition
$$
X = X_1 \sqcup \cdots \sqcup X_m.
$$
We write $x \sim x'$ if $x$ and $x'$ belong to the same $X_i$. 

Fix a field $k$. 
We associate to $X$ a product of matrix algebras 
$$
A(X)=A_1 \times \cdots \times A_m,
$$
each $A_i$ being 
the ring of $|X_i| \times |X_i|$ matrices over $k$ with rows and columns indexed by the elements of 
$X_i$. Equivalently, $A(X)$ consists of all functions $f:X \times X \to k$ such that $f(x,x')=0$
if $x \not\sim x'$ with multiplication 
$$
(fg)(x,x'') = \sum_{x' \in X} f(x,x')g(x',x'').
$$
In more suggestive notation, writing $f_{xx'}$ for $f(x,x')$, the product is
$$
(fg)_{x,x''} = \sum_{x' \in X} f_{x,x'}g_{x',x''}.
$$
Equivalently, $A(X)$ is the groupoid algebra associated to $X$ where the objects are the elements of $X$ 
and there is a single arrow $x \to x'$ whenever $x \sim x'$.

\begin{prop}
\label{prop.psets.mm.algs}
Let $X = X_1 \sqcup \cdots \sqcup X_m$ and $Y= Y_1 \sqcup \cdots \sqcup Y_n$ be finite partitioned
sets. Let $\tau:X \to Y$ be a function and define 
$$
X^j_i:=X_i \cap \tau^{-1}(Y_j)
$$ 
for all $i,j$. 
Suppose the restriction $\tau:X^j_i \to Y_j$ is bijective whenever $X^j_i$ is non-empty.  
Then the linear map  $\theta:A(Y) \to A(X)$ defined by  
$$
\theta(f)_{x,x'}:= \begin{cases}
				f_{\tau x,\tau x'} & \text{if $x \sim x'$}
				\\
				0 &  \text{if $x \not \sim x'$}
			\end{cases}
$$			
for all $f \in A(Y)$ and $x,x' \in X$ is an algebra homomorphism.
\end{prop}
\begin{pf}
First, $\theta(f)$ belongs to $A(X)$ because $ \theta(f)_{x,x'}=0$ if $x \not \sim x'$. If $x,z \in X$, then
\begin{align*}
\theta(fg)_{xz} & =\sum_{y \in Y} f_{\tau x,y} g_{y,\tau z}
\\
\big(\theta(f)\theta(g)\big)_{xz} & = \sum_{\substack{w \in X \\ x \sim w \sim z} } \theta(f)_{xw} \theta(g)_{wz} 
=\sum_{\substack{w \in X \\ x \sim w \sim z} }f_{\tau x,\tau w} g_{\tau w,\tau z}.
\end{align*}
To show  $\theta$ is an algebra homomorphism we must show the two sums on the right are
the same for all $x,z \in X$. If $x \not\sim z$, then $\theta(fg)_{xz}$ and $\big(\theta(f)\theta(g)\big)_{xz}$ are zero because the functions belong to $A(X)$. 
We therefore assume that $x \sim z$. Both right-hand sums are zero
if $\tau x \not\sim \tau z$ so we can, and do, assume $\tau x \sim \tau z$. Hence there exist $i$ and $j$ such that $x,z \in X_i$ and $\tau x,\tau z \in Y_j$. Thus $x,z \in X_i^j$. 

Comparing the two sums, one is over 
$$
\{w \in X \; | \; x \sim w \sim z\} \cap \{w \in X \; | \; \tau x \sim \tau w \sim \tau z\} =X_i \cap \tau^{-1}(Y_j) = X_i^j
$$
and the other is over $\{y \in Y \; | \; \tau x \sim y \sim \tau z\} = Y_j$. Since $X_i^j$ is not empty, $\tau:X^j_i \to Y_j$ is a bijection so the two sums are the same, whence
$$
\theta(fg) = \theta(f)\theta(g).
$$

Finally we check that $\theta$ sends the identity in $A(Y)$ to the identity in $A(X)$. The identity in $A(Y)$ is the 
function $e:Y \times Y \to k$ such that $e_{yy}=1$ for all $y \in Y$ and $e_{y,y'}=0$ if $y \ne y'$.
Hence 
$$
\theta(e)_{x,x'}:= \begin{cases}
				e_{\tau x,\tau x'} & \text{if $x \sim x'$}
				\\
				0 &  \text{if $x \not \sim x'$}
			\end{cases}
			\quad = \quad
			\begin{cases}
				1 & \text{if $x \sim x'$ and $\tau x = \tau x'$}
				\\
				0 &  \text{otherwise.}
			\end{cases}
$$			
If $x=x'$, then $x \sim x'$ and $\tau x = \tau x'$ so $\theta(e)_{xx} =1$. 
Suppose $x \ne x'$. If $x \not\sim x'$, then $\theta(e)_{xx'}=0$. If $x \sim x'$, then $\tau x \ne \tau x'$ because of the injectivity of the maps $\tau:X^j_i \to Y_j$ so $\theta(e)_{xx'}=0$.Hence $\theta(e)$ is the identity. 
\end{pf}

The next result is obvious.

\begin{prop}
\label{prop.Brat}
The Bratteli diagram \cite{B} associated to the map $\theta:A(Y) \to A(X)$  
consists of two horizontal rows with $n$ vertices on the top row  labelled by the numbers 
$|Y_1|, \ldots, |Y_n|$ and $m$ vertices on the lower row labelled by the numbers 
$|X_1|, \ldots, |X_m|$ and an edge from the top vertex $|Y_j|$ to the lower vertex $|X_i|$ whenever
$X_i^j \ne \varnothing$.
\end{prop}

\subsubsection{Example: Penrose sequences} 
\label{ssect.Connes.alg}
Let $X(n)$ be the set of Penrose sequences $z_0 \ldots z_n$, i.e., length-$(n+1)$ sequences
 of 0s and 1s in which the subsequence 11 does not appear. Let 
$\tau:X(n+1) \to X(n)$ be the map $\tau(z_0\ldots z_{n+1})=z_0 \ldots z_n$. Let $X(n)_0$ be the set of sequences in $X(n)$ ending in $0$, and  $X(n)_1$ the set of 
sequences ending in $1$. It is easy to verify that $\tau$ satisfies the condition in 
Proposition \ref{prop.psets.mm.algs}. Since
$|X(n+1)_0|=f_{n+2}$ and $|X(n+1)_1|= f_{n+1}$,
the Bratteli diagram is 
$$
  \UseComputerModernTips
\xymatrix{
f_{n+1}  \ar@{-}[dr] \ar@{-}[d] & f_n \ar@{-}[dl]
\\
f_{n+2} & f_{n+1}
}
$$
meaning that the associated algebra homomorphism
\begin{equation}
\label{eq.fib.system}
\theta:M_{f_{n+1}}(k) \times M_{f_{n}}(k) \to M_{f_{n+2}}(k) \times M_{f_{n+1}}(k) 
\end{equation}
 is
\begin{equation}
\label{eq.fib.maps}
\theta_n(a,b) = \Bigg( \begin{pmatrix} a & 0   \\ 0 & b \end{pmatrix}, a \Bigg).
\end{equation}

\subsubsection{Failure of functoriality} 

Given partitioned sets and maps $X  \stackrel{\tau}{\longrightarrow} Y  \stackrel{\sigma}{\longrightarrow} Z$ in which both $\tau$ and $\s$ have the property in Proposition \ref{prop.psets.mm.algs}, $\s\tau:X \to Z$ need not have the property in Proposition \ref{prop.psets.mm.algs}. An example of this is
$$
X=X_1=\{000,100,010,110\}, \quad Y_1=\{00,10\}, Y_2=\{10,11\}, \quad Z=Z_1=\{0,1\}
$$
with $\tau$ and $\s$ both being the map {\it drop the final digit}.

\subsection{(0,1)-sequences with at most $r$ consecutive 1s}
\label{sect.01s}
Fix $r\ge 1$. 
For all $n \ge 0$, and $0 \le i,j \le r$ we define
$$
\begin{array}{l}
{\sX}:=  \{\hbox{sequences $z_0z_1\ldots $  of 0s and 1s without $r+1$ consecutive 1s}\};
\\
{\sX}(n):=  \{\hbox{sequences $z_0z_1\ldots z_{n}$  of 0s and 1s without $r+1$ consecutive 1s}\};
\\
{\sX}(n)_i:=  \{\hbox{$z_0z_1\ldots z_{n} \in \sX(n)$  that end $01^i$}\} \quad \hbox{if $0\le i \le n$};
\\
{\sX}(n)_{n+1}:=  \{1\ldots 1= 1^{n+1}\};
\\
\tau:  \sX(n+1) \to \sX(n), \qquad \tau(z_0\ldots z_{n+1}):= z_0\ldots z_{n};
\\
X(n+1)_i^j:= \sX(n+1)_i \cap \tau^{-1}\big(\sX(n)_j\big).
\end{array}
$$

\begin{lemma}
\label{lem.conds.good}
The conditions in Proposition \ref{prop.psets.mm.algs} hold for $\tau:\sX(n+1) \to \sX(n)$
with respect to the partition 
$$
\sX(n) = \sX(n)_0 \sqcup \ldots \sqcup \sX(n)_r.
$$
\end{lemma}
\begin{pf}
The set $\sX(n+1)^j_0 = \{z_0\ldots z_{n-j-1}01^j 0 \in X(n+1)\}$ is non-empty if and only if $n\ge j$,
and in that case $\tau$ maps it bijectively to $\sX(n)_j$.

Suppose $i \ge 1$. Then $\sX(n+1)_i^j$ is non-empty if and only if $  i=j+1 \le n+1$ and in that case it
consists of all words  $z_0\ldots z_{n-i}01^i$ in $X(n+1)$ and $\tau$ maps it 
bijectively to $\sX(n)_{i-1}$. 
\end{pf}

It is convenient to  make the following definitions:
\begin{align*}
\ve:= &\hbox{the empty sequence of 0's and 1's},
\\
\sX(-1):=&\{\ve\},
 \\
\sX(-1)_{0}:=&\{\ve\},  
\\
\sX(n):=&\varnothing \qquad \hbox{for $n \le -2$}.
\end{align*}

 \begin{lemma}
 \label{lem.cn}
 Define $c_n$, $n \in \ZZ$, by 
\begin{itemize}
  \item 
  $c_i=0$ for $i \le -2$, 
  \item 
  $c_{-1}=1$,
  \item 
  $c_{n}=c_{n-1}+c_{n-2}+ \cdots +c_{n-r-1}$ for $n \ge 0$.
\end{itemize}  
 Thus $c_n=2^n$ for $0 \le i \le r$. 
 For all $n \in \ZZ$,  
\begin{equation}
\label{cn.card.Xn}
 |\sX(n)| = c_{n+1}.
\end{equation}
 \end{lemma}
 \begin{pf}
 It is clear that (\ref{cn.card.Xn}) holds for $n \le -1$.
 
The recurrence relation for $c_n$ implies $c_n=2^n$ for $0 \le n \le r$.
 Furthermore, if $-1 \le n \le r-1$, $\sX(n)$ consists of all length-$(n+1)$ sequences of 0s and 1s so
 $|\sX(n)|= 2^{n+1}$. Hence $|\sX(n)|=c_{n+1}$ for $-1 \le n \le r-1$. It is also clear that $|\sX(r)|=2^{r+1}-1=c_{r+1}$. Hence (\ref{cn.card.Xn}) holds for all $n \le r$. 
 
 We now assume $n \ge r$ and (\ref{cn.card.Xn}) is true for $n$.
 We will prove (\ref{cn.card.Xn}) holds with $n+1$ in place of $n$. 
 
 The sets $\sX(n)_i$ and $\sX(n)_i^j$ are only defined for $0 \le i,j \le r$ so the qualifications $j \le n$
 and $i \le n+1$ in the proof of Lemma \ref{lem.conds.good} are satisfied when $n \ge r$. 
 Because $n \ge r$, the  proof of Lemma \ref{lem.conds.good} therefore tells us that 
 \begin{enumerate}
  \item 
   $|\sX(n+1)_0^j| = |\sX(n)_j|$  for all $j$;
  \item 
   $|\sX(n+1)_i^j| = \varnothing$ if $i \ge 1$ and $j \ne i-1$;
  \item 
   $|\sX(n+1)_i^{i-1}| =  |\sX(n)_{i-1}|$ if $i \ge 1$.
\end{enumerate}

\underline{Claim:} If $n \ge r$ and $0 \le i \le r$, then 
\begin{equation}
\label{cni.card.Xni}
|\sX(n)_i|=c_{n-i}.
\end{equation}
\underline{Proof}: Since $\sX(r)_i$ consists of sequences $z=z_0\ldots z_{r-i-1}01^i$ where $z_0\ldots z_{r-i-1}$ is an arbitrary sequence of length $r-i$, $|\sX(r)_i|=2^{r-i}=c_{r-i}$. Hence  (\ref{cni.card.Xni}) 
holds for $n=r$ and all $i=0,\ldots,r$. 
Now assume that (\ref{cni.card.Xni}) holds for some $n \ge r$ and all $i$. We have  
\begin{align*}
|\sX(n+1)_0|& = |\sX(n+1)_0^0\sqcup \sX(n+1)_0^1 \sqcup  \cdots \sqcup \sX(n+1)_0^r|
\\
&= |\sX(n+1)_0^0|+  |\sX(n+1)_0^1|+ \cdots +  |\sX(n+1)_0^r|
\\
&= |\sX(n)_0|+  |\sX(n)_1|+ \cdots +  |\sX(n)_r|
\\
& = c_{n}+c_{n-1}+\cdots + c_{n-r}
\\
&=c_{n+1}.
\end{align*}
If $i \ge 1$, then
$$
|\sX(n+1)_i| = |\sX(n+1)^{i-1}_i|=|X(n)_{i-1}|=c_{n+1-i}.
$$
Hence (\ref{cni.card.Xni}) holds with $n+1$ in place of $n$ and for all $i=0,\ldots,r$. By induction, (\ref{cni.card.Xni}) holds for  all $n\ge r$ and all $i=0,\ldots,r$. $\lozenge$
 
It follows from the claim that
\begin{align*}
|\sX(n+1)|& = |\sX(n+1)_0\sqcup \sX(n+1)_1 \sqcup  \cdots \sqcup \sX(n+1)_r|
\\
& = c_{n+1}+c_{n}+\cdots + c_{n+1-r}
\\
&=c_{n+2}
\end{align*}
so (\ref{cn.card.Xn}) is true for all $n\ge 0$. 
 \end{pf}

\subsection{The algebras $A(n)$ and their direct limit}
For $n \ge 0$ we define
\begin{equation}
\label{defn.An}
A(n): =A(\sX(n-1)).
\end{equation}
Thus $A(0)=k$. 

The inverse limit of the directed system of sets $\cdots \to \sX(n+1) \to \sX(n) \to   \cdots$
is $\sX$ together with the maps $\rho_n:\sX \to \sX(n)$ given by $\rho_n(z_0z_1\ldots)=
z_0\ldots z_n$. Applying the procedure described in Proposition \ref{prop.psets.mm.algs} to the inverse system
$$
\sX =\limproj_n \sX(n) \longrightarrow \quad \cdots \stackrel{\tau}{\longrightarrow} \sX(2) \stackrel{\tau}{\longrightarrow} \sX(1) 
\stackrel{\tau} {\longrightarrow} \sX(0) \stackrel{\tau} {\longrightarrow} \sX(-1)
$$
produces a  directed system of $k$-algebras and $k$-algebra homomorphisms
\begin{equation}
\label{defn.Ans}
k= A(0)  \stackrel{\theta}{\longrightarrow} A(1) \stackrel{\theta}{\longrightarrow} A(2) 
\stackrel{\theta} {\longrightarrow} A(3) \stackrel{\theta} {\longrightarrow}  \cdots \quad 
\longrightarrow A=\liminj_n A(n).
\end{equation}
The following picture for $r=3$ serves as a mnemonic to understand the construction of 
each $A(n)$ and the 
maps $\tau$ that induce the homomorphisms $A(n) \to A(n+1)$:
$$
  \UseComputerModernTips
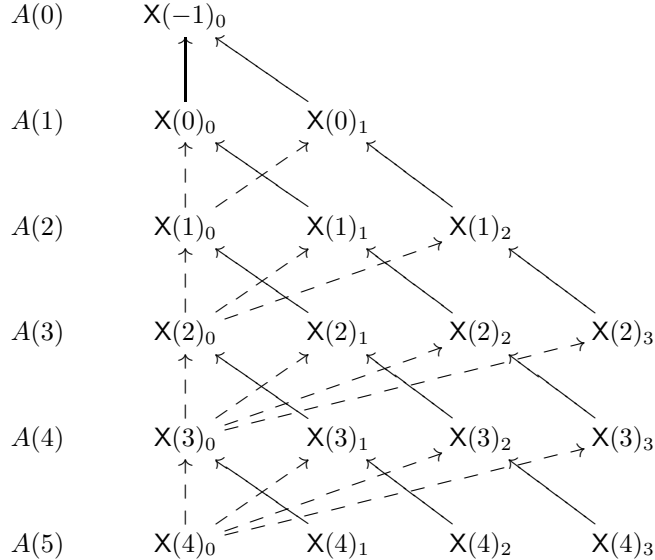
\xymatrix{
A(0) & \sX(-1)_0 
\\
A(1)  & \sX(0)_0 \ar[u] & \sX(0)_1 \ar[ul]
\\
A(2) & \sX(1)_0 \ar@{-->}[u]   \ar@{-->}[ur] & \sX(1)_1 \ar[ul]  & \sX(1)_2 \ar[ul]
\\
A(3) & \sX(2)_0 \ar@{-->}[u]   \ar@{-->}[ur]  \ar@{-->}[urr] & \sX(2)_1  \ar[ul] & \sX(2)_2 \ar[ul] &  \sX(2)_3 \ar[ul]
\\
A(4) & \sX(3)_0 \ar@{-->}[u]   \ar@{-->}[ur]  \ar@{-->}[urr] \ar@{-->}[urrr]  & \sX(3)_1 \ar[ul]  & \sX(3)_2 \ar[ul] & \sX(3)_3 \ar[ul] 
\\
A(5) & \sX(4)_0 \ar@{-->}[u]   \ar@{-->}[ur]  \ar@{-->}[urr] \ar@{-->}[urrr]   & \sX(4)_1 \ar[ul]  & \sX(4)_2 \ar[ul] & \sX(4)_3 \ar[ul] 
}
$$
The meaning of the solid arrows is that $\tau\big(\sX(n)_i\big) \subset \sX(n-1)_{i-1}$ when $i \ge 1$ (in fact $\tau$ restricts to a bijection between $\sX(n)_i$ and $\sX(n-1)_{i-1}$).
The dashed arrows from $\sX(n)_0$ to $\sX(n-1)_j$ 
mean that  $\tau\big(\sX(n)_0\big) \cap \sX(n-1)_j \ne \varnothing$
for all $j=0,\ldots, \min\{r,n\}$. 

We now reinterpret this diagram, more precisely its generalization for an arbitrary $r \ge 1$, to 
obtain a description of the homomorphisms. 

\begin{prop}
\label{prop.cn}
Let $c_n$, $n \in \ZZ$, be the numbers defined in Lemma \ref{lem.cn}. 
The Bratteli diagram for $A$ is repeated copies of 
$$
  \UseComputerModernTips
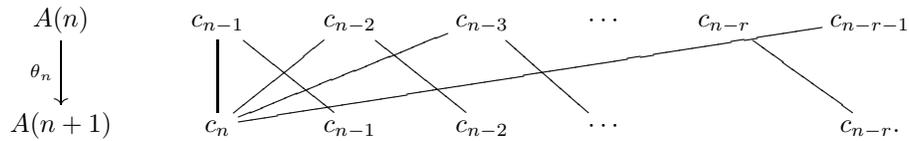
\xymatrix{
A(n) \ar[d]_{\theta_n} & c_{n-1}  \ar@{-}[dr]  \ar@{-}[d] & c_{n-2} \ar@{-}[dl] \ar@{-}[dr]  & c_{n-3}   \ar@{-}[dr]   \ar@{-}[dll]  & \cdots & c_{n-r} \ar@{-}[dr]  & c_{n-r-1}  \ar@{-}[dlllll] 
\\
A(n+1) & c_n & c_{n-1} &  c_{n-2}   & \cdots &   & c_{n-r}. 
}
$$
\end{prop}
\begin{pf}
Suppose $n \ge r$. 
Then $\sX(n-1)=\sX(n-1)_0 \sqcup \ldots  \sqcup \sX(n-1)_r$ so $A(n)=A(\sX(n-1))$ is a product of $r+1$ matrix algebras of sizes $|\sX(n-1)_0|, \ldots, |\sX(n-1)_r|$. Similarly, $A(n+1)$ is a product of $r+1$ matrix algebras of sizes $|\sX(n)_0|, \ldots, |\sX(n)_r|$. 
By Proposition \ref{prop.Brat}, the Bratteli diagram for the map $A(n) \to A(n+1)$ has an edge from the   top vertex labelled by $c_{n-j}=|\sX(n-1)_j|$ to the lower vertex labelled  
$c_{n-i}=|\sX(n)_i|$ whenever $\sX(n)_i^j \ne \varnothing$.

By the proof of Lemma \ref{lem.conds.good}, $\sX(n)^j_0$ is non-empty for all $j$ so there is an edge from the top vertex labelled $c_{n-j}$ to the bottom vertex labelled $c_n$  for all $j$.

 Suppose $i \ge 1$.  The proof of Lemma \ref{lem.conds.good} shows that $\sX(n)_i^{i-1}$ is 
 non-empty  and $\sX(n)_i^{j}$ is empty if $j\ne i-1$; i.e., there is an edge from the top vertex labelled $c_{n-i}=|\sX(n-1)_{i-1}|$ to the  lower vertex labelled  $c_{n-i}=|\sX(n)_i|$ and no other edges from the top row 
 to the  lower vertex labelled  $c_{n-i}=|\sX(n)_i|$. 

This completes the proof that the Bratteli diagram is as claimed for $n \ge r$. 

Now suppose $0 \le n \le r-1$. Then $A(n)$ is a product of $n+1$ matrix algebras of sizes 
$|\sX(n-1)_0|=c_{n-1},\ldots,|\sX(n-1)_{n}|=c_{-1}$; however, $c_i=0$ for $i \le -2$ so we can, and will, think of $A(n)$ as a product of $r+1$ matrix algebras of sizes $c_{n-1},\ldots,c_{-1}, c_{-2},\ldots,c_{n-r-1}$. With this convention, and a generalization of the diagram just prior to the statement of this proposition one sees the Bratteli diagram is as claimed for $0 \le n \le r-1$. 
\end{pf} 

\begin{thm}
The $k$-algebra $S$ defined in section \ref{sect.S} is isomorphic to $A= \liminj A(n)$.
\end{thm}
\begin{pf}
{\it If} $b_n=c_n$ for all $n \ge 0$,  the Bratteli diagram in Proposition \ref{prop.cn} for the directed system of $A(n)$s is the same as  the Bratteli diagram  for the directed system of $S_n$s in Proposition \ref{prop.Brat.Sn} and therefore $A \cong S$ by Elliott's Theorem \cite{E}. 

By Proposition \ref{prop.B}, the sequence $b_n=\dim_kB_n$, $n \ge 0$, satisfies the same recurrence relation as the sequence $c_n$, $n \ge 0$. Since $b_i=2^i=c_i$ for $0 \le i \le r$, $c_n=b_n$ for all $n \ge 0$. 
\end{pf}

The Bratteli diagrams for the cases $r=1$, $r=2$, and $r=3$, are
$$
  \UseComputerModernTips
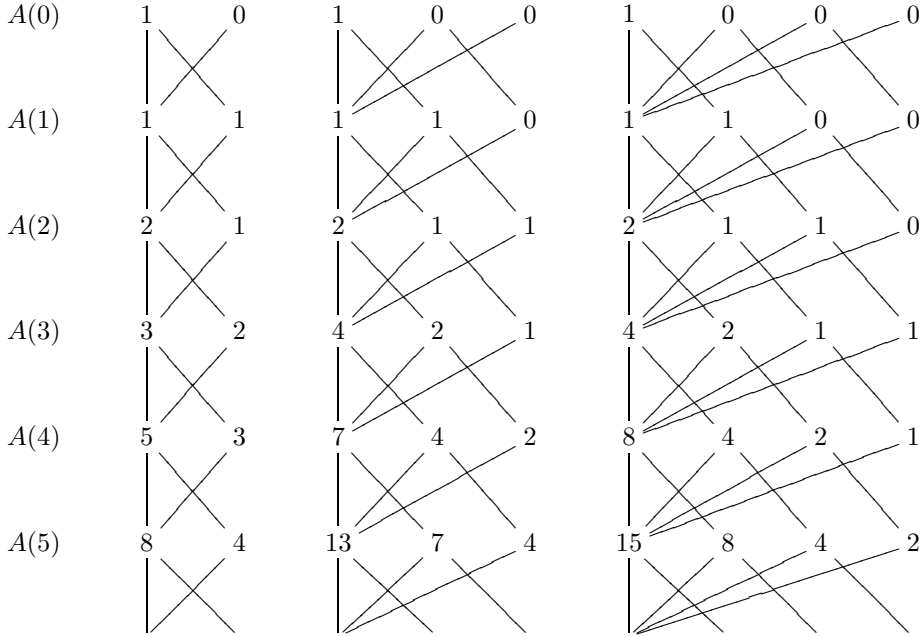
\xymatrix{
A(0) & 1  \ar@{-}[dr] \ar@{-}[d] & 0 \ar@{-}[dl]  & 
 1  \ar@{-}[dr]  \ar@{-}[d] & 0 \ar@{-}[dl] \ar@{-}[dr]  & 0  \ar@{-}[dll] &
 1  \ar@{-}[dr]  \ar@{-}[d] & 0 \ar@{-}[dl] \ar@{-}[dr]  & 0  \ar@{-}[dll]  \ar@{-}[dr] & 0  \ar@{-}[dlll] 
\\
A(1) & 1  \ar@{-}[dr] \ar@{-}[d] &1 \ar@{-}[dl] &
 1   \ar@{-}[dr]  \ar@{-}[d]  &1 \ar@{-}[dl] \ar@{-}[dr]  & 0   \ar@{-}[dll] &
 1  \ar@{-}[dr]  \ar@{-}[d] & 1 \ar@{-}[dl] \ar@{-}[dr]  & 0  \ar@{-}[dll] \ar@{-}[dr] & 0  \ar@{-}[dlll] 
\\
A(2) & 2  \ar@{-}[dr] \ar@{-}[d] &1 \ar@{-}[dl] &
 2  \ar@{-}[dr]  \ar@{-}[d]  & 1  \ar@{-}[dl] \ar@{-}[dr] &1  \ar@{-}[dll] &
 2  \ar@{-}[dr]  \ar@{-}[d] & 1 \ar@{-}[dl] \ar@{-}[dr]  & 1  \ar@{-}[dll] \ar@{-}[dr] & 0  \ar@{-}[dlll] 
\\
A(3) & 3  \ar@{-}[dr] \ar@{-}[d] &2 \ar@{-}[dl] &
 4    \ar@{-}[dr]   \ar@{-}[d]  & 2 \ar@{-}[dl] \ar@{-}[dr]  &1  \ar@{-}[dll] &
 4  \ar@{-}[dr]  \ar@{-}[d] & 2 \ar@{-}[dl] \ar@{-}[dr]  & 1  \ar@{-}[dll] \ar@{-}[dr] & 1  \ar@{-}[dlll] 
\\
A(4) & 5  \ar@{-}[dr] \ar@{-}[d] &3 \ar@{-}[dl] &
 7    \ar@{-}[dr]   \ar@{-}[d]  & 4 \ar@{-}[dl] \ar@{-}[dr]  &2  \ar@{-}[dll] &
 8  \ar@{-}[dr]  \ar@{-}[d] & 4 \ar@{-}[dl] \ar@{-}[dr]  & 2  \ar@{-}[dll]  \ar@{-}[dr]  & 1 \ar@{-}[dlll] 
\\
A(5) & 8  \ar@{-}[dr] \ar@{-}[d] &4 \ar@{-}[dl] &
 13    \ar@{-}[dr] \ar@{-}[d] & 7 \ar@{-}[dl] \ar@{-}[dr]   &4  \ar@{-}[dll] &
 15  \ar@{-}[dr]  \ar@{-}[d] & 8 \ar@{-}[dl] \ar@{-}[dr]  & 4  \ar@{-}[dll]  \ar@{-}[dr]  & 2 \ar@{-}[dlll] 
\\
&& &&&&& & &&&& &
}
$$
The Fibonacci sequence appears in the left-most column of the left-most Bratteli diagram and 
the tribonacci sequence occurs in the left-most column of the middle diagram. 

\subsection{Infinite paths in the Bratteli diagram}

Fix $r$ and consider the corresponding Bratteli diagram. For example, set $r=3$ and consider the right-most of the three diagrams above. Relabel the vertices by placing 0 in each left-most vertex and label all other vertices with 1. Replace each edge by a downward pointing arrow. An infinite path in the Bratteli  diagram is a sequence $z_0z_1\ldots$ in which $z_i$ is a vertex at level $i$ and there is a
downward arrow $z_i \to z_{i+1}$ (the $0^{\th}$ level is the top one). The set of all infinite paths is the set of all infinite sequences of 0s and 1s such that $r+1$ consecutive 1s never occur. The infinite paths in the Bratteli diagram are in natural bijection with the points in $\sX$.

\section{The Grothendieck group of $\qgr B$}
\label{sect.K0}
 
As before, $B=k\langle x,y\rangle/(y^{r+1})$.

\subsection{}
Let $\a$ be the smallest positive real root of the irreducible polynomial $t^{r+1}-t^r-\cdots -t-1$
(we say more about $\a$ in section \ref{sect.pisot}). 
We will prove $K_0(\qgr(B))$ is isomorphic to $\ZZ[\a]$ as an ordered abelian group with order unit $[\cO]$
corresponding to $1 \in \ZZ[\a]$, $[\cO(-1)]$ corresponding to $\a$, and $[\cF(-1)]=\a[\cF]$ for all $\cF \in \qgr(kQ)$. We will also describe the structure of $K_0(\qgr(B))$ as a $\ZZ[t^{\pm 1}]$-module. 

The $\ZZ[t^{\pm 1}]$-module structure arises as follows.

Given a group $G$ acting as automorphisms of an abelian category $\sA$, $K_0(\sA)$ becomes a $\ZZ G$-module by $g\cdot [\cF]:=[g \cdot \cF]$.  Because $B$ is a $\ZZ$-graded ring the degree shift functor $M \rightsquigarrow M(1)$ gives $K_0(\gr B)$ the structure of a $\ZZ[t,t^{-1}]$-module with $t$ acting as $(-1)$. 
The Serre twist, $\cF \rightsquigarrow \cF(1)$, on $\qgr B$  gives $K_0(\qgr B)$ the structure of a 
$\ZZ[t,t^{-1}]$-module with $t$ acting as $(-1)$.  

\subsection{}
\label{sect.kQ.0}
  
  One might try to compute the Grothendieck group of the category $\qgr(B)$ by using the localization sequence $K_0(\fdim B) \to K_0(\gr B) \to K_0(\qgr B) \to 0$. However, the global dimension of $B$ is 
  infinite so we cannot conclude that the natural map $K_0(\sP(\gr B)) \to K_0(\gr B)$ induced by the inclusion of the exact subcategory $\sP(\gr B)$ of projective objects in $\gr B$ is an isomorphism. It therefore seems a difficult task to compute  $K_0(\gr B)$.
  
However, there is an inclusion of $B$ in the path algebra of a quiver $Q$, $f:B \to kQ$,
with the property that $kQ \otimes_B -$ induces an equivalence of categories $\qgr(B) \to \qgr(kQ)$. Since $kQ$ has finite global dimension we can compute $K_0(\gr (kQ))$ easily then use the localization sequence $K_0(\fdim(kQ)) \to K_0(\gr (kQ)) \to K_0(\qgr (kQ)) \to 0$ to compute $K_0(\qgr (kQ))$ and hence $K_0(\qgr B)$.

The methods and results in sections \ref{sect.kQ.1} and \ref{sect.kQ.2} parallel those in \cite{HS1} 
and \cite{HS2}. The results in those papers relate a monomial algebra to the path algebra of 
its Ufnarovskii graph but here $Q$ is not the Ufnarovskii graph of $B$.

\subsection{}
\label{sect.kQ.1}
Fix an integer $r \ge 1$ and  let $Q$ be the quiver
$$
 \UseComputerModernTips
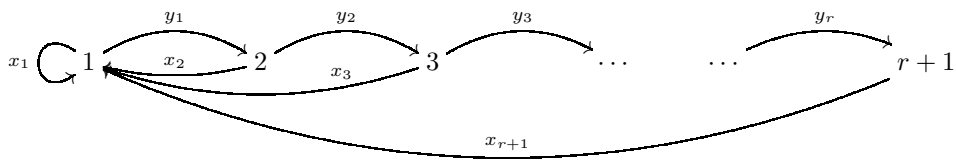
\xymatrix{
1 \ar@(ul,dl)[]_{x_1} \ar@/^1pc/[rr]^{y_1} && 2 \ar@/^.4pc/[ll]_{x_2} \ar@/^1pc/[rr]^{y_2} && 3   \ar@/^1pc/[llll]_<<<<<<{x_3}  \ar@/^1pc/[rr]^{y_3} && \cdots & \cdots \ar@/^1pc/[rr]^{y_r} && r+1  \ar@/^3pc/[lllllllll]_{x_{r+1}}
}
$$
\vskip .3in
\noindent
We write $kQ$ for the path algebra of $Q$.  We write $Q_0$ for the set of vertices in $Q$ and $Q_1$ for the set of arrows. If $a$ is an arrow we write $s(a)$ for its start and $h(a)$ for its head: visually,
$s(a) \stackrel{a}{\longrightarrow} h(a)$.

\subsubsection{Notation.}
Let $p$ and $q$ be paths in $Q$ such that $p$ ends where $q$ starts. We adopt the convention that $pq$ denotes the path {\it traverse $p$ then $q$}. We write $e_i$ for the trivial path, i.e., the primitive idempotent, at vertex $i$.

\begin{prop}
\label{prop.map.f}
There is an injective algebra homomorphism $f:B \to kQ$ such that
\begin{align*}
f(x) & = x_1+ \cdots + x_{r+1}
\\
f(y) & = y_1+\cdots +y_r.
\end{align*}
\end{prop}
\begin{pf}
Since $(y_1+\cdots +y_r)^{r+1}=0$, $f$ {\it is} a well-defined algebra homomorphism.

A  basis for $B$ is provided by the set
$$
W:=\{\hbox{words in $x$ and $y$ not containing $y^{r+1}$ as a subword}\}.
$$
The function $\psi:\{\hbox{paths in $Q$ that begin at vertex 1} \} \to W$ defined by 
$$
\psi(p):=\begin{cases}
	 & \text{the word  obtained by removing the subscripts}
	 \\
	 & \text{from the arrows $x_i$ and $y_j$ that appear in $p$}.
	 \end{cases}
$$
is a bijection.  

To show $f$ is injective it suffices to show that $\{f(w_\l) \; | \; \l \in \L\}$ is linearly independent 
whenever $\{w_\l \; | \; \l \in \L\}$ is a finite subset of $W$. Indeed, it is enough to show that $\{e_1f(w_\l) \; | \; \l \in \L\}$ is linearly independent. We will prove this by showing that $\{e_1f(w_\l) \; | \; \l \in \L\}$
consists of different paths.  

We will do this by showing that 
\begin{enumerate}
  \item[(a)]
  $e_1f(w)$ is a path (necessarily beginning at vertex 1) and
  \item[(b)] 
  $\psi(e_1f(w))=w$ for all $w \in W$. 
\end{enumerate} 
(The linear independence of $\{e_1f(w_\l) \; | \; \l \in \L\}$ follows from (b).)
We argue by induction on $\deg w$. It is clear that $e_1f(1)=1$ and $\psi(e_1f(1))=1$ so (a) and (b)
are true when $\deg w=0$. 

Suppose (a) and (b) are true for $w \in W$. 

Then $e_1f(wx)=e_1f(w)f(x)=e_1f(w)(x_1+\cdots+x_{r+1})$; the fact that $e_1f(w)$ is a path   implies $e_1f(w)f(x)=e_1f(w)x_j$ for some $j$ . Hence (a) holds for $wx$. Because  $\psi(e_1f(w))=w$, 
$\psi(e_1f(w)x_j)=\psi(e_1f(w))x=wx$, i.e., (b) holds for $wx$. 

Now consider $e_1f(wy)=e_1f(w)(y_1+\cdots+y_r)$. If the path $e_1f(w)$ ends at a vertex 
$i \ne r+1$, then $e_1f(wy)= e_1f(w)y_i \ne 0$. 
The path $e_1f(w)$  begins at vertex 1 so if it were to end at $r+1$ it would equal 
$p'y_1\ldots y_r$ for some path $p'$ which would imply
$w=\psi(e_1f(w))=\psi(p'y_1\ldots y_r)=\psi(p')y^r$ whence $wy=\psi(p')y^{r+1}=0$, i.e., $wy \notin W$.
We conclude that  $e_1f(w)f(y)$ is non-zero and therefore equal to $e_1f(w)y_i$ for some $i$. Hence (a) holds for $wy$. So does (b) because $\psi(e_1f(wy))=\psi(e_1f(w)y_i)=\psi(e_1f(w))y=wy$.

This completes the proof that (a) and (b) hold for all $w \in W$.  Item (b) implies that $\{e_1f(w_\l) \; | \; \l \in \L\}$ consists of different paths. The injectivity of $f$ follows. 
\end{pf}

\subsection{}
\label{sect.kQ.2}
We will follow the strategy in \cite{HS1} and \cite{HS2} and use the following theorem of Artin and Zhang to show that the map $f:B \to kQ$ induces an equivalence of categories $\QGr(B) \to \QGr(kQ)$
via the functor $kQ \otimes_B -$.

 \begin{prop}
\label{prop.AZ2}
\cite[Prop 2.5]{AZ}
\label{prop.AZ}
Let $B$ and $D$ be $\NN$-graded $k$-algebras such that 
$\dim_k B_i < \infty$ and $\dim_k D_i < \infty$
for all $i$. Let $f:B \to D$ be a homomorphism of graded $k$-algebras. 
If $\ker f$ and $\coker f$ belong to $\Fdim B$, then $D \otimes_B -$ induces an equivalence of categories 
\begin{equation}
\label{eq.AZ}
\QGr B \to \frac{\Gr D}{\sT_B}
\end{equation}
where
$$
\sT_{B}=\{M\in \Gr D\;|\; M_B\in \Fdim B\}.
$$
\end{prop}

We now  proceed to prove that the hypotheses of Proposition \ref{prop.AZ} hold for the map 
$f:B \to kQ$ and that
the category $\sT_B$ in (\ref{eq.AZ}) is equal to $\Fdim(kQ)$.

\begin{lem}
cf. \cite[Lemma 2.2]{HS2}.
Let $f:B \to kQ$ be the homomorphism in Proposition \ref{prop.map.f} and write $(kQ)B_n$ for the 
left ideal of $kQ$ generated by $f(B_n)$.  For all $n \ge 0$,
$$
(kQ)B_n=kQ_{\ge n}.
$$
\end{lem}
\begin{pf}
Write $D=kQ$. 

We will prove that $D_0B_n=D_n$ for all $n \ge 0$. This is certainly true for $n=0$.

Since $e_if(x)=x_i$ and $e_if(y)=y_i$, $D_0B_1=D_1$. 
We now argue by induction on $n$. If $D_0B_{n-1}=D_{n-1}$, then 
$$
D_0B_n = D_0(B_1)^n=D_0(B_1)^{n-1}B_1=D_{n-1}B_1=D_{n-1}D_0B_1=D_{n-1}D_1=D_n.
$$
This completes the proof that $D_0B_n=D_n$ for all $n \ge 0$. Hence $DB_n=D_{\geq n}$.
\end{pf}

\begin{lem}
\label{lem.coker}
The cokernel of the homomorphism $f:B \to kQ$ in Proposition \ref{prop.map.f} belongs to $\Fdim B$. 
\end{lem}
\begin{pf}
Write $D=kQ$. 

By \cite[Lemma 2.2]{HS1}, $\coker(f)$ belongs to $\Fdim(B)$ if $f(B_{m})D_0\subset f( B_m)$ and 
$ f(B_m) D_1 \subset f(B_{m+1})$ for some $m \in \NN$. We will show that 
$f(B_{r+1})D_0\subset f( B_{r+1})$ and $ f(B_{r+1}) D_1 \subset f(B_{r+2})$.

\underline{Claim:} 
If the letter $x$ appears in $w \in W$, there is an arrow $a$ such that $f(w)=da$ for some $d \in kQ$. 

\underline{Proof:} 
If $w=w''xw'$ for $w',w'' \in W$, then $f(w)=f(w'')(x_1+\cdots+x_{r+1})f(w') = f(w''x)e_1f(w')$ because all arrows $x_i$ end at vertex 1. However, $e_1f(w')$ is a path by the statement (a) 
in the proof of Proposition \ref{prop.map.f}. In particular $e_1f(w')=d'a$ for some arrow $a$ and some path 
$d'$. Hence $f(w)=f(w''x)d'a$ for some arrow $a$.
$\lozenge$

Suppose $w \in B_{r+1}$. 
Since $y^{r+1}=0$, at least one $x$ appears in $w$. Therefore $f(w)=da$ for some arrow $a$ and some $d \in kQ$. 

Hence $f(w)e_i=dae_i$ is either 0 or $da$. In either case,  
$f(w)e_i \in f(B_{r+1})$. Therefore $f(B_{r+1})D_0\subset f( B_{r+1})$.

Now consider $f(w)D_1$.  
If $a$ ends at vertex $r+1$, then $f(w)Q_1 = \{0,dax_{r+1}\}$ and $dax_{r+1}=daf(x)=f(wx) \in 
f(B_{r+2})$.
If $a$ ends at vertex $i \ne r+1$, then $f(w)Q_1 = \{0,dax_i, day_i\}$ and $dax_i=f(wx) 
\in f(B_{r+2})$ and $vay_i=f(wy) \in f(B_{r+2})$. This completes the proof that 
$f(B_{r+1})D_1\subset f( B_{r+2})$. 

Therefore \cite[Lemma 2.2]{HS1} implies $\coker(f)$ belongs to $\Fdim(B)$
\end{pf}

The proof of the next result is essentially identical to that of \cite[Prop.2.3]{HS2}.

\begin{prop}
\label{prop.whew}
$\sT_B=\Fdim (kQ)$.
\end{prop}

The next theorem now follows from Propositions \ref{prop.map.f}, \ref{prop.AZ}, \ref{prop.whew}, and 
Lemma \ref{lem.coker}.
 
\begin{thm}
Let $f:B \to kQ$ be the homomorphism in Proposition \ref{prop.map.f}. The functor $kQ \otimes_B-$
induces an equivalence of categories 
$$
F:\QGr(B) \stackrel{\equiv}{\longrightarrow} \QGr(kQ).
$$
Furthermore, 
\begin{enumerate}
  \item 
 $F$ restricts to an equivalence between $\qgr(B)$ and $\qgr(kQ)$; 
  \item 
  $F(\cF(1)] \cong F(\cF)(1)$ for all $\cF \in \QGr(B)$;
  \item 
  if $\cO$ is used to denote both the image of ${}_BB$ in $\QGr(B)$ and  ${}_{kQ}kQ$ in $\QGr(kQ)$,
  then $F(\cO)=\cO$.  
\end{enumerate}
\end{thm}

\subsection{The Pisot number $\a$}
\label{sect.pisot}
Let $\a$ be an algebraic integer belonging to $\RR$. We call $\a$ a {\sf Pisot number} if $\a >1$ and
all its Galois conjugates have modulus $<1$.  We refer the reader to 
\cite{FS} for basic information about Pisot numbers.

By \cite{Brauer}, the polynomial $f(x):=x^{r+1}-x^r-\cdots -x -1$ is irreducible for all 
$r \ge 1$ and 
\begin{equation}
\label{defn.alpha}
\a:= \hbox{the unique positive real root of $f(x)$}
\end{equation}
is a Pisot number. 
As remarked in \cite{BMS}, the $\a$s are monotonically increasing and converge to 2 as $r$ increases.

In \cite[p.3]{BMS}, the $(r+1)^{st}$ multinacci number $\omega_{r+1}$ is defined to be the positive solution of the equation
$$
1 -t -t^2 - \cdots -t^{r+1}.
$$
Thus, $\omega_{r+1}=\a^{-1}$. 

When $r=1$, $\a:=\frac{1+\sqrt{5}}{2}$, the Golden Ratio.
When $r=2$, $\a$ is the Tribonacci constant, the limit of the ratio
$t_n/t_{n-1}$ of successive terms in the Tribonacci sequence defined by $t_n=t_{n-1}+t_{n-2}+t_{n-3}$
and $t_0=1$, $t_1=2$, $t_2=4$. Recall, $t_n$ is the dimension of the degree $n$ component of the algebra $k\langle x,y\rangle/(y^3)$. 
 
The following result, which concerns the growth of the 
 coefficients of the power series expansion of the rational function $H_B(t)$, is probably implicit, if not explicit, in the literature.
 
 \begin{prop}
 \label{prop.growth.B}
 Let $B=k\langle x,y\rangle/(y^{r+1})$. Then 
 $$
 \lim_{n \to \infty} \frac{\dim B_n}{ \a^n} = \frac{\a^{r+2}}{1+2\a+3\a^2+\cdots+(r+1)\a^r}.
 $$  
 \end{prop}
 \begin{pf}
 By Proposition \ref{prop.B},
 $$
 H_B(t)=\frac{1+t+\cdots+t^r}{1-t-\cdots-t^{r+1}}.
 $$
 The denominator is $t^{r+1}f(t^{-1})$ which has a simple zero at $\a^{-1}$ so
 $$
 H_B(t)=\frac{c\a}{1-\a t} +\frac{h(t)}{g(t)}
 $$
 for some $c \in \RR$, where $(1-\a t)g(t)$ is the  denominator of $H_B(t)$ and $h(t)$ is some polynomial. It follows that 
 $$
c \a = \lim_{t \to \a^{-1}} H_B(t)(1-\a t).
 $$
Since $H_B(t)(1-\a t)$ is a  ratio of polynomials that vanish at $\a^{-1}$, l'H\^opital's rule tells us that
$$
c \a \,  = \lim_{t \to \a^{-1}} H_B(t)(1-\a t) =  \frac{\a^{r+2}}{1+2\a+3\a^2+\cdots+(r+1)\a^r}.
$$
On the other hand, 
$$
H_B(t)=c\a \sum_{i=0}^\infty \a^it^i + \sum_{i=0}^\infty d_it^i
$$
where the second term is the formal power series expansion of $h(t)/g(t)$. However, since $\a$ is a Pisot number all the zeroes of $g(t)$ lie outside the unit disk and 
therefore $\lim_{i \to \infty} d_i =0$.  But
$$
\frac{\dim B_n}{\a^n}=  c\a  + \frac{d_n}{\a^n}
$$
so the result follows. 
 \end{pf}

\subsection{}
We are now ready to compute the Grothendieck group $K_0(\qgr(B))$, first as a $\ZZ[t^{\pm 1}]$-module
then as an ordered abelian group with order unit.

\begin{thm}
\label{thm.K0.1}
Let $Q$ be the quiver at the beginning of section \ref{sect.kQ.1}. There is
an isomorphism
$$
K_0(\qgr(B)) \cong K_0(\qgr(kQ)) \cong \frac{\ZZ[t^{\pm 1}]}{(1-t-\cdots - t^{r+1})}
$$
of $\ZZ[t^{\pm 1}]$-modules under which $[\cO] \mapsto 1$ and $[\cF(-1)]=t[\cF]$ for all $\cF \in \qgr(B)$.
\end{thm}
\begin{pf}
We will use the localization sequence $K_0(\fdim(kQ)) \to K_0(\gr (kQ)) \to K_0(\qgr (kQ)) \to 0$ to compute $K_0(\qgr (kQ))$.

Let $\sP$ be the exact subcategory of projective objects in $\gr(kQ)$. The inclusion $\sP \to \gr(kQ)$
induces a group homomorphism $K_0(\sP) \to K_0(\gr(kQ))$. 
Since $kQ$ has finite global dimension (=1) this map is an isomorphism. 
The indecomposable projectives in $\gr(kQ)$ are the graded modules $P_i:=(kQ)e_i$, $i \in Q_0$,  
and their twists $P_i(n)$, $n \in \ZZ$.  Let $p_i$ denote the class $[P_i]$ in $K_0(\gr(kQ))$. 
Thus $K_0(\gr(kQ))$ is a free $\ZZ[t^{\pm 1}]$-module with basis
$p_1,\ldots,p_{r+1}$.

Let $V_i$, $i \in Q_0$, be the 1-dimensional simple  graded left
$kQ$-module supported at $i$ and concentrated in degree 0. The $V_i$s belong to 
$\fdim(kQ)$. Every object in $\fdim(kQ)$ has a composition series whose factors are 
isomorphic to twists of the $V_i$s. By d\'evissage $K_0(\fdim(kQ))$ is the free $\ZZ[t^{\pm 1}]$-module with basis the classes $[V_i]$, $1 \le i \le r+1$. Let $v_i$ be the image of $[V_i]$
in $K_0(\gr(kQ))$ under the natural map $K_0(\fdim(kQ)) \to K_0(\gr (kQ)$. The localization sequence 
tells us that $K_0(\qgr(kQ))$ is isomorphic as a $\ZZ[t^{\pm 1}]$-module to 
$K_0(\gr (kQ))$ modulo the submodule generated by $v_1,\ldots,v_{r+1}$.

 The minimal projective resolution of $V_i$ in $\gr(kQ)$ is 
$$
0 \longrightarrow \bigoplus_{a \in h^{-1}(i)} P_{s(a)}(-1) \stackrel{(\cdot a)}{\longrightarrow} P_i \longrightarrow V_i \longrightarrow 0
$$
so  
$$
v_i=[P_i] - \sum_{a \in h^{-1}(i)} [P_{s(a)}(-1)] = p_i - t\sum_{a \in h^{-1}(i)} p_{s(a)}.
$$

Hence $K_0(\qgr(kQ))$ is isomorphic to $\oplus_{i=1}^{r+1} \ZZ[t^{\pm 1}]p_i$ modulo the 
the relations
$$
p_i = t\sum_{a \in h^{-1}(i)} p_{s(a)}, \qquad 1 \le i \le r+1.
$$
These relations are $p_1=t(p_1+\cdots + p_{r+1})$ and $p_{i}=tp_{i-1}$ for $i=2,\ldots,r+1$. 
The relations imply that $p_i=t^{i-1}p_1$ for $i=2,\ldots,r+1$ and $p_1=t(1+t+\cdots+t^r)p_1$. 

Hence $K_0(\qgr(kQ))$ is isomorphic to the rank one free module $\ZZ[t^{\pm 1}]p_1$ modulo 
the relation $p_1=t(1+t+\cdots +t^r)p_1$. 
Since $kQ = P_1 \oplus \cdots \oplus P_{r+1}$, $[\cO]=(1+ t +\cdots + t^r)p_1=t^{-1}p_1$ so $K_0(\qgr(kQ))$ is also generated by $[\cO]$ as a $\ZZ[t^{\pm 1}]$-module. Hence 
$K_0(\qgr(kQ))$ is isomorphic to the rank one free module $\ZZ[t^{\pm 1}]\cdot[\cO]$ modulo 
the relation $[\cO]=t(1+t+\cdots +t^r)[\cO]$.  By sending $[\cO]$ to 1 we obtain an isomorphism
$$
K_0(\qgr(kQ)) \cong \frac{\ZZ[t^{\pm 1}]}{(1-t-\cdots - t^{r+1})}
$$
of $\ZZ[t^{\pm 1}]$-modules. It is  a tautology that $[\cF(-1)]=t[\cF]$.
\end{pf}

We need a special case of the following result of Frougny and Solomyak.

\begin{thm}
\cite[Thm. 2]{FS}
\label{thm.FS}
Let $f(x)=x^m - a_1 x^{m-1}-\cdots -a_m$ where $a_i \in \ZZ$ and $a_1 \ge a_2 \ge \cdots \ge a_m>0$.
Then 
\begin{enumerate}
  \item 
  $f(x)$ has a unique positive root, $\b$;
  \item 
   $\b$ is a Pisot number;
  \item{}
  $$
  \ZZ[\b^{-1}]\cap \RR_{\ge 0}  = \bigcup_{s=0}^\infty \Big\{\sum_{i=-s}^s a_i \b^i \; \Big\vert \; a_i \in \{0,1,\ldots,[\b]\}
\Big\}.
  $$
\end{enumerate}
\end{thm}

To describe the order structure on $K_0(\qgr(B)$ we embed it in $(\RR,+)$.

\begin{thm}
\label{thm.K0}
Let $\a$ be the unique positive  root of $x^{r+1}-x^r-\cdots -x -1$. 
\begin{enumerate}
\item{}
  There is an injective group homomorphism $K_0(\qgr B) \to \RR$ given by 
  $[\cO] \mapsto 1$ and $[\cO(-n)] \mapsto \a^n$ that maps $K_0(\qgr B)$ 
  isomorphically to $\ZZ[\a]$.
  \item{}
  The image of the positive cone in $K_0(\qgr B)$ under the homomorphism 
  in (1) is $\RR^+ \cap \ZZ[\a] = \ZZ[\a^{-1}] \cap \RR^+$.  
\end{enumerate}
\end{thm}
\begin{pf}
(1)
Since the constant term in the minimal polynomial of $\a$ is 1,  $\ZZ[\a]=\ZZ[\a^{-1}]$. 
The minimal polynomial of $\a^{-1}$ is $1-t-\cdots - t^{r+1}$ so there is a $\ZZ$-algebra isomorphism
$$
 \frac{\ZZ[t^{\pm 1}]}{(1-t-\cdots - t^{r+1})} \to \ZZ[\a], \qquad t \mapsto \a^{-1}.
$$
Combining this with the description of $K_0(\qgr B)$ in Theorem \ref{thm.K0.1}
gives a $\ZZ[t^{\pm 1}]$-module isomorphism $K_0(\qgr B) \to \ZZ[\a]$ where $t$ acts on $\ZZ[\a]$
as multiplication by $\a^{-1}$.

(2)
If $a_i \in \NN$, then $a_i\a^i \leftrightarrow [\cO(-i)^{\oplus a_i}]$ under the
isomorphism in (1) so Theorem \ref{thm.FS} implies that given $\c \in \ZZ[\a]\cap \RR_{\ge 0}$  
there is an object  $\cF \in \qgr B$ such that  the isomorphism in (1) sends $[\cF]$ to $\c$.

It remains to show that if $\cF \in \qgr B$, then $[\cF] \in \ZZ[\a] \cap \RR^+$.\footnote{At present we
are unable to prove this without using the equivalence $\qgr B \equiv \mod S$ 
but we hope to find an alternative proof
that relies only on arguments involving $B$ and $kQ$.}
Since $\qgr B$ is equivalent to $\mod S$, $K_0(\qgr B) \cong K_0(S)=\liminj K_0(S_n)$. 
We will write elements of $\ZZ^{r+1}$ as column vectors.

From now on we consider the directed system $K_0(S_r) \to K_0(S_{r+1}) \to \cdots$.
Hence in all that follows $n$ is assumed to be $\ge r$.

The algebra $S_n$ is a product of $(r+1)$ matrix algebras of sizes $b_n,b_{n-1},\ldots,b_{n-r}$ so
we make the identification
$$
K_0(S_n) =  \ZZ^{r+1}
$$
with $K_0(S_n)_{\ge 0} = \NN^{r+1}$ and order unit $[S_n]=(b_n,\ldots,b_{n-r})^{\sf T}$.

It follows from the Bratteli diagram for the inclusion $\phi_n:S_n \to S_{n+1}$ in Proposition \ref{prop.Brat.Sn} that
the induced map $K_0(S_n) \to K_0(S_{n+1})$ is left multiplication by 
\begin{equation}
\label{M.matrix}
M:= \begin{pmatrix}
  1&1&1& \cdots  & & 1   \\
    1  & 0 & 0 &   \cdots  &  & 0 \\
   0  & 1 & 0 &   \cdots    && 0 \\
   \vdots & &&&& \vdots \\
     0  & 0 & 0 &   \cdots & 1  & 0 \\
\end{pmatrix}.
\end{equation}
(One can check that $M\cdot [S_n]=M (b_n,\ldots,b_{n-r})^{\sf T} =[S_{n+1}]$  because 
$b_{n+1}=b_n+\cdots+b_{n-r}$.)   
The directed system  $K_0(S_r) \to K_0(S_{r+1}) \to \cdots$ is therefore isomorphic to 
$$
\ZZ^{r+1} \stackrel{M}{\longrightarrow}  \ZZ^{r+1} \stackrel{M}{\longrightarrow} \ZZ^{r+1} \stackrel{M}{\longrightarrow} \cdots.
$$
The vector $\vul=(v_0,\ldots,v_r)$ where
\begin{align*}
v_r & = \a^r,\\
 v_{r-1}& =\a^r+\a^{r-1}, \\
 \vdots \phantom{1}  & \phantom{=\;} \phantom{=\;}  \vdots \\
 v_1& =\a^r+\a^{r-1}+\cdots+\a, \\ 
 v_0&=\a^{r+1},
 \end{align*}
is a left $\a$-eigenvector for $M$. Since $\vul M=\a\vul$, the diagram
\begin{equation}
\label{dlim.K0}
 \UseComputerModernTips
\xymatrix{
\ZZ^{r+1} \ar[d]_{\vul} \ar[r]^M & \ZZ^{r+1} \ar[d]_{\a^{-1}\vul} \ar[r]^M & \ZZ^{r+1} \ar[d]_{\a^{-2}\vul} \ar[r]^M & \cdots
\\
\RR \ar@{=}[r] & \RR \ar@{=}[r] &\RR \ar@{=}[r] & \cdots
}
\end{equation}
commutes. 
Since the entries of $\vul$ are positive, the vertical maps $\a^{-n}\vul\cdot$
send $\NN^{r+1}$ to $\RR_{\ge 0}$. Hence $K_0(S)_{\ge 0} \subset \RR_{\ge 0}$. Since the entries of 
$\a^{-n}\vul$ belong to $\ZZ[\a]$, the vertical maps
send $\NN^{r+1}$ to $\ZZ[\a]$. Hence $K_0(S)_{\ge 0} \subset \RR_{\ge 0} \cap \ZZ[\a]$.
\end{pf}
 
\subsubsection{Remark}
The order unit in $K_0(S)=\ZZ[\a]$ is the image of $[S_r] =(2^r, \ldots,2,1)^{\sf T}$ under the left-most vertical map in (\ref{dlim.K0}). One can replace $\vul$ by a non-zero scalar multiple of itself in the vertical maps of (\ref{dlim.K0}) to make the order unit 1.

\subsection{}
Let $\cM \in \qgr B$. There is an exact sequence 
$$
0 \to \bigoplus_i \cO(-i)^{c_i} \to \bigoplus_i \cO(-i)^{d_i} \to \cM \to 0 
$$
for a suitable set of non-negative integers $c_i$ and $d_i$ only finitely many of which are non-zero. 
We define 
$$
p_\cM(t):=\sum_i(d_i-c_i)t^i.
$$
Under the isomorphism $K_0(\qgr B) \to \ZZ[\a]$, $[\cO] \mapsto 1$, $[\cO(-1)] \mapsto \a$, the image of $[\cM]$ in  $\ZZ[\a]$ is $p_{\cM}(\a^{-1})$. 

\begin{cor}
If $\cM$ and $\cN$ belong to $\qgr B$ and $p_{\cM}(\a^{-1}) \le p_{\cN}(\a^{-1})$, then 
$\cN \cong \cM \oplus \cF$ for some $\cF \in \qgr B$.
\end{cor}
\begin{pf}
The hypothesis implies that  $[\cN]-[\cM] \ge 0$ so Theorem \ref{thm.K0} implies that $[\cN]-[\cM]=[\cF]$ for some $\cF \in \qgr B$. Hence $[\cM \oplus \cF] =[\cN]$ and therefore 
$\cM \oplus \cF \oplus \cG \cong \cN \oplus \cG$ for some $\cG \in \qgr G$ by the definition of $K_0$.
But multi-matrix algebras are unit regular rings so $S$, being a direct limit of such, is
also unit-regular. By \cite[Thm. 4.5]{vNR} we can cancel $\cG$ to obtain the claimed result.
\end{pf}

\section{The perspective of non-commutative algebraic geometry}

In the language of non-commutative algebraic geometry, $\Projnc B$, which is defined implicitly
by declaring that $\Qcoh(\Projnc B):=\QGr B$, is an {\it affine} nc-scheme 
because its category of ``quasi-coherent sheaves'' is equivalent to  a category of modules over a ring.

\subsection{Skyscraper sheaves associated to elements of $\sX$}

As in section \ref{sect.01s},  we write ${\sX}$ for the set of sequences $z_0z_1\ldots $  of 0s and 1s having at most $r$ consecutive 1s. Elements $z$ and $z'$ in $\sX$ are equivalent,  denoted  $z \sim z'$, if 
$z_{\ge n}=z'_{\ge n}$ for some $n$.

 Given $z \in {\sf X}$, let $M_z$ be the point module defined in section 
\ref{ssect.2nd.thm} and let $\cO_z$ be $M_z$ viewed as an object in $\QGr B$; i.e., $\cO_z = \b^* M_z$.

\begin{prop}
\label{prop.sky}
Let $z,z' \in {\sf X}$. Then $\cO_z \cong \cO_{z'}$ if and only if $z \sim z'$.
\end{prop}
\begin{pf}
See  section \ref{ssect.2nd.thm}.
\end{pf}

\begin{prop}
\label{prop.F2}
Suppose $k=\FF_2$. If $M$ is a point module for $B$, then $\b^*M \cong \cO_z$ for some $z \in \sX$.
\end{prop}
\begin{pf}
Let $e_0,e_1,\ldots$ be a homogeneous basis for $M$ with $\deg e_i =i$. Define $z_i \in \FF_2$ by
$y.e_i=z_ie_{i+1}$. Then $z_0z_1z_2 \ldots$ is a sequence of 0s and 1s and every 1 is followed by a 0 because $y^2=0$.
\end{pf}

Proposition \ref{prop.sky}, sets up a bijection between the equivalence classes of the Penrose tilings
and certain ``points'' on $\Projnc B$.

\begin{prop}
If $k \not\cong \FF_2$ there is a point module $M$, such that $\b^*M$ is not isomorphic
to any $\cO_z$, $z \in {\bf P}$.  
\end{prop}
\begin{pf}
Let $a \in k-\{0,1\}$. Let $\l=\l_0\l_1\l_2 \ldots$ be a sequence in which each $\l_i$ belongs to $\{0,1,a\}$
with the restriction that if $\l_i \ne 0$, then $\l_{i+1}=0$. Let $M$ be the graded vector space with basis 
$e_0,e_1,\ldots$ and $\deg e_i =i$. Make $M$ a graded left $B$-module by declaring that
$$
x.e_i=e_{i+1} \qquad y.e_i:=\l_i e_{i+1}.
$$
If $a$ appears infinitely often in $\l$, then $\b^*M$ is not isomorphic to any $\cO_z$.
\end{pf}

\subsection{$\Projnc B$ is a ``smooth non-commutative curve''}
\label{ssect.curve}

Let $X$ be an irreducible algebraic variety. Then $X$ is a smooth 
algebraic curve if and only if $\Ext^2_{\Qcoh X}$ is identically zero and $\Ext^1_{\Qcoh X}$ is not.

Given that fact, the next result suggests that  $\Projnc B$ might be thought of as a 
smooth non-commutative curve, albeit of a  strange kind.

\begin{prop}
$\Ext^2_{\QGr B}$ is identically zero and $\Ext^1_{\QGr B}$ is not.
\end{prop}
\begin{pf}
Since $S$ is a direct limit of semisimple $k$-algebras, the main result in \cite{Bers} tells us that the global homological dimension of $S$ is at most 1, i.e., $\Ext^2_{S}$ is identically zero. Since $S$ is not a semisimple ring it does not have global dimension 0. 
\end{pf}

\subsection{``Dense points'' in $\Projnc B$}
\label{ssect.dense}

The space of Penrose tilings modulo isometry is a coarse topological space, i.e., every point in it is dense.
 
In section \ref{ssect.top.P} (see Proposition \ref{prop.links}), we observed that
 $\Ext^1_{\QGr B}(\cO_z,\cO_{z'} ) \ne 0$ for all $z,z' \in {\sX}$. 
 We now explain why this can be interpreted as a ``density'' result.

First, we think of the isomorphism classes of the simple objects in $\QGr B$ as playing the role of ``points'' on the non-commutative curve $\Projnc B$. 

Let $P$ be a finite topological space satisfying the $T_0$ condition. The category of such topological 
spaces is equivalent to the category of finite posets. 
 We make  $P$ a poset by declaring that $x \le y$ if  $x \in \overline{\{y\}}$. 
  
 Let $I(P)$ be the incidence algebra of $P$. 
 The isomorphism classes of simple left $I(P)$-modules are in bijection with the points of $P$. 
 If $V_x$ and $V_y$ are simple left modules associated to distinct points 
 $x,y \in P$, then $\Ext^1_{I(P)}(V_x,V_y) \ne 0$ if and only if  $y$ {\it covers} $x$ which means that
 $x < y$ and if $x \le z \le y$, then $z \in \{x,y\}$. (Ladkani's result  \cite[Lemma 2.7]{Lad} that the category of sheaves of finite dimensional vector spaces on $P$ is equivalent to the category of finite dimensional right $I(P)$-modules persuades the author that $I(P)$ is a ``coordinate ring'' of the topological space $P$.)
 
  In particular, 
 $$
 \Ext^1_{I(P)}(V_x,V_y) \ne 0 \qquad \Longrightarrow \qquad x \in \overline{\{y\}}.
 $$
  Conversely, if $x \in \overline{\{y\}}$, there is a sequence of points $x=x_0<x_1< \cdots < x_n =y$ such that $\Ext^1_{I(P)}(V_{x_i},V_{x_{i+1}}) \ne 0$ for all $i$.
  
Since $\cO_z$ depends only on the equivalence class $[z]$ of $z$ in ${\sX}/\!\!\sim$, this result about incidence algebras suggests we might interpret the fact that 
 $\Ext^1_{\QGr B}(\cO_z,\cO_{z'} ) \ne 0$ for all $z,z' \in {\sX}$ as saying that 
 ``every $[z]$ is in the closure of every $\{[z']\}$'', i.e., every point in 
${\sX}/\!\!\sim$ is dense. As we said, this is indeed the case.

\appsectionA{Penrose tilings}

We have used Gr\"unbaum and Shepard's book \cite{GS} as our standard reference for tilings. We have also found the books by Senechal \cite{Sen} and Sadun \cite{Sad} helpful.  This appendix follows \cite[Ch. 10]{GS}
quite closely.

\smallskip

\noindent
{\bf A.1. Tilings of the plane.}
For our purposes the following definitions suffice.
A {\sf tiling} $\cT=\{T_1,T_2,\ldots\}$ of the plane is a countable family of closed subsets $T_i \subset \RR^2$ such that
\begin{enumerate}
  \item 
  the union of the $T_i$s is $\RR^2$;
   \item
  each $T_i$ is a polygon, not necessarily convex, and is the closure of its interior;
  \item 
  the interiors of the $T_i$s are pairwise disjoint;
  \item{}
  if $T_i \cap T_j \ne \varnothing$, then  $T_i \cap T_j$ is an edge of $T_i$ and an edge of $T_j$.
\end{enumerate}  
Let $\cT=\{T_1,T_2,\ldots\}$ and $\cT'=\{T'_1,T'_2,\ldots\}$ be tilings of  the plane.
We say $\cT_1$ and $\cT_2$ are {\sf congruent} if 
 there is an isometry $\s$ of the plane such that $\{\s T_1,\s T_2,\ldots \}=\{T'_1,T'_2,\ldots\}$.
 We say $\cT_1$ and $\cT_2$ are {\sf equal} or {\sf the same} if there is a positive real number $\l$ such that
$\{\l T_1,\l T_2,\ldots \}$ is congruent to $\cT'$.

An isometry of the plane is called a {\sf symmetry} of $\cT$ if $\{\s T_1,\s T_2,\ldots \}=\{T_1,T_2,\ldots\}$.

\smallskip

 \noindent
{\bf A.2. Prototiles.}
A set of closed polygons  $p_1,\ldots,p_n$ is called a {\sf set of prototiles} for the tiling 
$\cT=\{T_1,T_2,\ldots\}$ if every $T_i$ is congruent, directly or by reflection, to some $p_j$. 

\smallskip

\noindent
{\bf A.3. Kites and darts.}
Let $\tau=\frac{1+\sqrt{5}}{2}$ and $\theta=\frac{\pi}{5}$. A quadrangle with sides of length $1,1,\tau,\tau$ and interior angles $2\theta$, $2\theta$, $2\theta$, $4\theta$ is called a {\sf kite}. 
A quadrangle with sides of length $1,1,\tau,\tau$ and interior angles $\theta$, $\theta$, $2\theta$, $6\theta$ is called a {\sf dart}. We color the vertices of the kite and dart as follows: on the kite the vertex at the angle 
$4\theta$ and the vertex opposite it are colored black and the other two vertices are colored white;
on the dart the vertex at the angle $6\theta$ and the vertex opposite it are colored white and the other two vertices are colored black. See the picture on \cite[p.539]{GS}.

A tiling $\cT$ of the plane is called a {\sf Penrose tiling}  if
\begin{enumerate}
  \item 
   the kite and dart form a set of prototiles for $\cT$;
  \item 
  whenever a point is a vertex of more than one tile  it is the same color on each of those tiles.
\end{enumerate}
 Condition (2) is called a {\sf matching rule}.
 
 \smallskip
 
 \noindent
{\bf A.4. $A$-tiles.}
   Cutting a kite along its axis of symmetry produces two
  congruent isoceles triangles that have sides of length $1,\tau,\tau$ and acute  interior angles. 
  Cutting a dart along its axis of symmetry produces two congruent isoceles triangles having sides of length $1,1,\tau$ and an obtuse interior angle. The area of a triangle coming from a kite is bigger than that of a triangle coming from a dart. 
  We call the two triangles {\sf $A$-tiles} and label the larger one with the letter $A$ and the smaller with the letter $a$. The vertices of the $A$-tiles are colored by inheriting the colors of the vertices on the kite and dart.
  
  Given a Penrose tiling $\cT$ we can cut each kite and dart in it along its axis to produce a new
  tiling of the plane by  triangles. The triangles $A$ and $a$ form a set of prototiles for the new tiling. 
  We call the new tiling an {\sf $A$-tiling}. 
  
We can reverse this process provided we impose a matching rule for the $A$-tiles that involves matching the color of common vertices {\it and} matching certain orientations of the sides of the triangles  \cite[pp.538-540]{GS}. A tiling $\cT'$ for which the two $A$-tiles form a set of 
prototiles is called an {\sf $A$-tiling} if the matching rules are satisfied. 
Such an $A$-tiling $\cT'$ can be converted to a Penrose tiling 
by kites and darts if we delete all edges that join vertices of the same color or, equivalently, by ``amalgamating
the $A$-tiles in mirror image pairs'' \cite[p.540]{GS}.

\smallskip

 \noindent
{\bf A.5. $B$-tiles.}  
 We can form an isoceles triangle with sides of length $\tau,\tau,1+\tau$ by placing an $A$-triangle next to 
 an $a$-triangle in such a way that coincident vertices have the same color. We label this new triangle with
  the letter 
 $B$. The color of the vertices of $B$ and the orientation of its sides are inherited from those of $A$ and $a$. 
 
 Given an $A$-tiling $\cT$ we can delete certain edges so as to 
 amalgamate adjacent $A$- and $a$-tiles in a unique way 
 to form $B$-triangles. All the $a$-tiles disappear in this process so we obtain
 a new tiling having $B$ and $A$ as prototiles. We now replace 
 the label $A$ by the letter $b$ so the new tiling has $B$ and $b$ as prototiles. 
 The triangles $B$ and $b$ are called {\sf $B$-tiles} and we call the new tiling a 
 {\sf $B$-tiling}.  The area of  $B$ is bigger than that of $b$. 
 
This process of passing from an $A$-tiling to a $B$-tiling is called {\sf composition}.
 
 A similar composition process may be applied to a $B$-tiling to produce an $A$-tiling: one deletes edges so as to amalgamate adjacent  $B$- and $b$-tiles, subject to matching rules, to form a triangle congruent to $\tau A$ by which we mean a copy of $A$ scaled up by a factor of $\tau$.  
 Every $b$ tile disappears in this process and one now has a tiling with prototiles $\tau A$ and $B=\tau a$
 We shrink the $\tau A$-tiling by a factor of $\tau$ to get a genuine $A$-tiling. If we wish to
 we can then delete all edges that join vertices of the same color to produce a genuine 
 Penrose tiling by kites and darts.

\smallskip

 \noindent
{\bf A.6. Penrose sequences.}
Following \cite[p.568]{GS}, we now explain how to assign to each Penrose tiling a sequence $z_0z_1\ldots$ of $0$s and $1$s.

Let $\cT$ be a Penrose tiling and $p \in \RR^2$ a point that does not lie on the edge or axis of symmetry of any  kite or dart. Cut each kite and dart along its axis of symmetry to form an $A$-tiling.  We set $z_0=0$ if $p$ is in an $A$-tile and set $z_0=1$ if it is in an $a$-tile. We now compose to obtain a $B$-tiling
and set $z_1= 0$ if $p$ is in a $B$-tile and set $z_1=1$ if it is in a $b$-tile. We now apply the composition process to this $B$-tiling to get a $\tau A$ tiling, shrink by a factor of $\tau$, and determine $z_2$ in the same way we determined $z_0$. Notice that $z_1=0$ if $z_0=1$ because all $a$-tiles are absorbed into $B$-tiles
when composing an $A$-tiling to get a $B$-tiling. Likewise, $z_2=0$ if $z_1=1$ because in composing a $B$-tiling to get an $A$-tiling all $b$-tiles are absorbed into $\tau A$-tiles. 

We repeat this ad infinitum to produce a  sequence $z_0z_1\ldots$ in which each 1 is followed by a 0.
Every infinite sequence of $0$s and $1$s without adjacent 1s is obtained from some Penrose tiling $\cT$ and 
point $p$ \cite[Thm. 10.5.9]{GS}. To prove this one must construct a Penrose tiling from 
each such sequence; the proceedure is explained at \cite[p.568]{GS}. Some sequences lead to 
a tiling that fills up only a half-plane or a sector with angle $\pi/5$;  one then fills the whole plane by reflection(s). 
 
The sequence depends on the choice of $p$. However, if one picks a different point $p'$, then after a certain number of compositions $p$ and $p'$ end up in the same tile whence the sequence for $(\cT,p)$ is eventually the same as the sequence for $(\cT,p')$. In the vernacular, the sequences have the same {\sf tail}. 
Conversely, two sequences that have the same tail correspond to two points in the same tiling. 

If $\s$ is a symmetry of $\cT$ the sequences associated to $(\cT,p)$ and $(\cT,\s p)$ are the same. Conversely, if the sequences for $(\cT,p)$ and $(\cT',p')$ are identical, then $\cT$ and $\cT'$ are the same and $p'$ is the image of $p$ under some symmetry of $\cT$. 

\appsectionB{Classification of all point modules for $\CC\langle x,y \rangle/(y^2)$}
 
  \noindent
{\bf B.1. Point module sequences.}
A {\sf point-module sequence} for $\CC\langle x,y\rangle/(y^2)$  is a sequence of points ${\bf a}=a_0,a_1,\ldots$ on the Riemann sphere with the property that if $a_i \ne \infty$, then $a_{i+1}=\infty$. Two such sequences ${\bf a}$ and ${\bf b}$ are {\sf equivalent} if there is an integer $n$ such that $a_i=b_i$ for all  $i \ge n$.

\smallskip

  \noindent
{\bf B.2. Point modules.}
If $M=\CC e_0 \oplus \CC e_1 \oplus \cdots$ is a point module for $\CC\langle x,y\rangle/(y^2)$ with $\deg e_i=i$ we associate to $M$ the point-module sequence ${\bf a}(M)=(a_0,a_1,\ldots)$ defined by 
$$
a_i =\begin{cases}
	 \infty & \text{if $ye_i=0$}
	 \\
	 xe_i=a_iye_i &  \text{if $ye_i\ne 0$.}
	 \end{cases}
$$
Conversely, if ${\bf b}$ is a point-module sequence there is a unique-up-to-isomorphism point module
$M$ such that ${\bf a}(M)={\bf b}$. 

\begin{prop}
If $M$ and $N$ are point modules, then $\b^*M \cong \b^* N$ if and only if ${\bf a}(M)$ is equivalent to 
${\bf a}(N)$.
\end{prop}

  \noindent
{\bf B.3. Semi-infinite paths.}
There is a map $\Psi$ from the set of point-module sequences to the set of semi-infinite paths 
$p=p_0p_1\ldots$ in the directed graph
$$
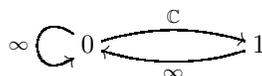
\UseComputerModernTips
\xymatrix{
0    \ar@/^/[rr]^{\CC}  \ar@(ul,dl)[]_\infty  && 1 \ar@/^/[ll]^\infty 
}
$$
that begin at the vertex labelled $0$.
If ${\bf a}=a_0a_1\ldots$ is a point-module sequence, define $\Psi({\bf a})=p_0p_1\ldots$ as follows:
if $a_i \ne \infty$, then $p_i$ is the arrow labelled $\CC$; there is then a unique way to assign arrows
labelled $\infty$ to the other $a_j$s so that $p_0p_1\ldots$ is an infinite path starting at $0$.

The fiber  of $\Psi$ over a path  $p_0p_1\ldots$ beginning at $0$ is a product of copies of 
$\CC$ and $\{\infty\}$. For example, the fiber that contains $\infty 1 \infty 2 \infty \infty 3 \infty \ldots$ is 
$$
\{\infty\} \times \CC \times  \infty  \times \CC \times \{\infty\} \times  \{\infty\}  \times \CC \times  \{\infty\}  \ldots
$$

 \end{document}